\documentclass[12pt,english]{article}
\usepackage[T1]{fontenc}
\usepackage[latin9]{inputenc}
\usepackage[letterpaper]{geometry}
\usepackage{color}
\usepackage{textcomp}
\usepackage{mathrsfs}
\usepackage{amsmath}
\usepackage{amsthm}
\usepackage{amssymb}

\makeatletter
\numberwithin{equation}{section}
\numberwithin{figure}{section}
\theoremstyle{plain}
\newtheorem{thm}{\protect\theoremname}[section]
\theoremstyle{definition}
\newtheorem{example}[thm]{\protect\examplename}
\theoremstyle{definition}
\newtheorem{defn}[thm]{\protect\definitionname}
\theoremstyle{plain}
\newtheorem{prop}[thm]{\protect\propositionname}
\theoremstyle{plain}
\newtheorem{cor}[thm]{\protect\corollaryname}
\theoremstyle{plain}
\newtheorem{lem}[thm]{\protect\lemmaname}

\@ifundefined{date}{}{\date{}}
\makeatother

\usepackage{babel}
\providecommand{\corollaryname}{Corollary}
\providecommand{\definitionname}{Definition}
\providecommand{\examplename}{Example}
\providecommand{\lemmaname}{Lemma}
\providecommand{\propositionname}{Proposition}
\providecommand{\theoremname}{Theorem}

\begin{document}
\title{A {\Large{}Constructive Real Projective Plane}}
\author{Mark Mandelkern}

\maketitle

\begin{quotation}
\noindent \textbf{\small{}Abstract.}{\small{} The classical theory
of plane projective geometry is examined constructively, using both
synthetic and analytic methods. The topics include Desargues's Theorem,
harmonic conjugates, projectivities, involutions, conics, Pascal's
Theorem, poles and polars. The axioms used for the synthetic treatment
are constructive versions of the traditional axioms. The analytic
construction is used to verify the consistency of the axioms; it is
based on the usual model in three-dimensional Euclidean space, using
only constructive properties of the real numbers. The methods of strict
constructivism, following principles put forward by Errett Bishop,
reveal the hidden constructive content of a portion of classical geometry.
A number of open problems remain for future studies. }\\
\end{quotation}
\noindent {\small{}2010}\emph{\small{} Mathematics Subject Classification.}{\small{}
Primary 51A05; Secondary 03F65.}{\small\par}

\noindent \emph{\small{}Key words and phrases.}{\small{} Projective
geometry, harmonic conjugates, projectivities, Pascal's Theorem, constructive
mathematics.}{\small\par}

\section*{{\large{}Introduction} \label{INTRODUCTION}}

In various forms\emph{,} the constructivist program goes back to Leopold
Kronecker (1823-1891), Henri Poincaré (1854-1912), L. E. J. Brouwer
(1881-1966) {[}Bro08{]}, and many others. The most significant recent
work, using the strictest methods, is due to Errett Bishop (1928-1983).
A large portion of analysis has been constructivized by Bishop in
\textit{Foundations of Constructive Analysis} {[}B67{]}; this treatise
also serves as a guide for constructive work in other fields. Expositions
of constructivist ideas and methods may be found in {[}B67, BB85{]};
see also {[}Sto70, R82, M85{]}. 

The initial phase of this program involves the rebuilding of classical
theories, using only constructive methods; the entire body of classical
mathematics is viewed as a wellspring of theories waiting to be constructivized. 
\begin{quotation}
Every theorem proved with {[}nonconstructive{]} methods presents a
challenge: to find a constructive version, and to give it a constructive
proof.

\noindent ~~~~~~~~~~~~~~~~~~~~~~~~~~~~~~~~~~~~~~~~~~~~~~~~~~~~~~~~~-
Errett Bishop {[}B67, p. x{]} 
\end{quotation}
\quad{}The present work is based on the classical theory of the real
projective plane. The classical theory is highly nonconstructive;
it relies heavily, at nearly every turn, on the\textit{ Law of Excluded
Middle}\emph{.} For example, it is classically assumed that a given
point is either on a given line, or not on the line \textemdash 	although
no finite routine is available for making such a determination; a
constructive development must utilize only the finite routines that
are specified in the axioms. 

\textcolor{black}{\emph{Background; classical.}}\textcolor{black}{{}
}Guides to the classical theory that have proven useful include works
by O. Veblen, J. W. Young, H. S. M. Coxeter, E. Artin, and G. Pickert
{[}VY10, Cox55, Art57, Pic75{]}. For a concise historical review,
with thorough references, see Cremona's preface {[}Cre73, pp. v-xii{]}.
An entertaining history of the subject will be found in {\small{}Lehmer's}
last chapter {[}Leh17, pp. 122-143{]}{\small{}. }{\small\par}

\emph{Background; constructive.} A. Heyting {[}H28, D90{]} has developed
a portion of the theory, based on axioms for projective space. A plane
being thus embedded in a higher dimension, this permits a proof of
Desargues's Theorem, and aids the verification of the properties of
harmonic conjugates. D. van Dalen {[}D96{]} has studied alternative
axioms involving the basic relations. 

Here we proceed in a different direction; we utilize axioms only for
a plane. Since there exist non-Desarguesian projective planes,\footnote{See, for example,  {[}Wei07{]}.}
this means that Desargues's Theorem must be taken as an axiom; it
will be used to establish the converse and the essential properties
of harmonic conjugates. The theory is developed further, as far as
conic sections, Pascal's Theorem, and the theory of polarity. We make
full use of duality in establishing some of the fundamental results. 

Adhering closely to the methods of strict constructivism, as introduced
by Errett Bishop {[}B67{]}, we eschew additional assumptions, such
as those of formal-logic-based intuitionism or recursive function
theory. For a full account of the distinctions between these various
types of constructivism, see {[}BR87{]}. 

\emph{Background; other work in constructive geometry. }For the constructive
extension of an affine plane to a projective plane, see {[}H59, D63,
M13a, M14{]}. For the constructive coördinatization of a plane, see
{[}M07{]}. 

The constructive geometry of M. Beeson {[}Bee10{]} uses Markov's Principle,\footnote{Markov's Principle and other nonconstructive principles will be discussed
in section 12.} which is accepted in recursive function theory, but not in the Bishop-type
strict constructivism that is followed in the present paper. M. Lombard
and R. Vesley {[}LomVes98{]} construct axioms for classical and intuitionistic
plane geometry, using recursive function theory. 

The work of J. von Plato {[}Pla95, Pla98, Pla10{]}, proceeding within
formal logic, is related to type theory, computer implementation,
and combinatorial analysis. For constructive plane hyperbolic geometry,
also within formal logic, see V. Pambuccian {[}Pam01{]}. The Bishop-type
constructive mathematics of the present paper works from a position
well-nigh opposite that of formal logic; for further comments on this
distinction, see {[}B65, B67, B73, B75{]}. 

\emph{Synthetic and analytic approaches. }We determine the constructive
possibilities of synthetic methods, using no axioms of order, in constructing
a real projective plane \emph{$\mathbb{P}$}. The analytic model $\mathbb{P}^{2}(\mathbb{R})$,
constructed in Euclidean space $\mathbb{R}^{3}$, is used to prove
the consistency of the axiom system. 

\emph{Axioms. }In choosing axioms for the projective plane $\mathbb{P}$,
we claim to adopt no new axioms, using only constructive versions
of the usual classical axioms. These axioms are all constructively
valid on the plane $\mathbb{P}^{2}(\mathbb{R})$, taking note of Bishop's
thesis, ``All mathematics should have numerical meaning'' {[}B67,
p. ix{]}\emph{. }No axioms of order are involved here; the constructive
investigation of ordered projective planes must be left for future
studies.\footnote{For a survey of classical axiomatic ordered geometry, see {[}Pam11{]}.
For a constructive theory of ordered affine geometry, see {[}Pla98{]}.} 

\emph{Bishop-type constructivism. }We follow the constructivist principles
put forward by Errett Bishop in 1967. This variety of constructivism
does not form a separate branch of mathematics, nor is it a branch
of logic; it is intended as an enhanced approach for all of mathematics.
For the distinctive characteristics of Bishop-type constructivism,
as opposed to intuitionism or recursive function theory, see {[}BR87{]}. 

\textit{Logical setting.} This work uses informal intuitionistic logic;
it does not operate within a formal logical system. For the origins
of modern constructivism, and the disengagement of mathematics from
formal logic, see Bishop's ``Constructivist Manifesto'' {[}B67,
BB85; Chapter 1{]}. Concerning the source of misunderstandings in
the mathematical community as to the methods and philosophy of constructivism,
see {[}B65{]}. 

We use intuitionistic logic only so far as to eschew use of the \emph{Law
of Excluded Middle} and its corollaries. Intuitionism, in the stronger
sense of Brouwer, introduces additional principles which are classically
false. In the opposite direction, recursive function theory limits
consideration to a restricted class of objects.\footnote{\noindent For more information concerning these alternative approaches
to constructivism, see {[}BR87{]}. } Constructive mathematics as proposed by Bishop leads down neither
of these extreme pathways. No special logical assumptions are made.
Avoiding the \emph{Law of Excluded Middle, }constructive mathematics
is a generalization of classical mathematics, just as group theory,
a generalization of abelian group theory, avoids the commutative law.
Every result and proof obtained constructively is also classically
valid. 

\emph{Results.} A fair portion of classical projective geometry is
found to have a solid constructive content, provided that appropriate
axioms, definitions and methods are used. 

It is necessary to avoid the classically ubiquitous method of treating
separately elements that are, or are not, distinct or separated from
other elements; constructively, elements typically arise lacking such
information. Thus, harmonic conjugates must be given a single definition
for all points on a line, without distinguishing between the base
points and other points; projectivities must be shown to exist for
any two given ranges, not knowing whether they are identical or distinct;
the polar of a point with respect to a conic must be defined uniformly
for any point on the plane, without treating points on the conic as
special cases. These requirements often reduce the role of the quadrangle
in definitions. 

Basing the theory only on axioms for a plane, with no axioms of space,
means that Desargues's Theorem must be taken as an axiom; the converse
is proved as a consequence. The necessity of ensuring that triangles
claimed to be perspective have the required distinctness is paid due
attention. Similar situations arise in proving the validity of the
harmonic conjugate construction, and the other main concepts. 

Once the basic definitions and properties are established constructively,
the theory proceeds rather smoothly, revisiting results obtained over
the centuries \textemdash{} now with constructive methods. 

\part{{\large{}Synthetic constructions} \label{PART. Synthetic}}

\noindent From a set of constructively phrased axioms, we derive the
basic properties of a real projective plane, including harmonic conjugates,
projectivities, conics, Pascal's Theorem, poles, and polars. 

\section{{\large{}Constructive methods}\label{SECTION 1. cxtv methods}}

\noindent One characteristic feature of the constructivist program
is meticulous use of the connective ``or\emph{''}. To prove ``$A\textnormal{ or }B$''
constructively, it is required that either we prove $A$, or we prove
$B$; it is not sufficient to prove the contrapositive $\neg(\neg A\text{ and }\neg B).$
For an illustration of this in geometry, the Brouwerian counterexample
given below will show that the statement ``Either the point $P$
lies on the line $l$, or $P$ lies outside $l$'' is constructively
invalid.

\emph{Constructively invalid statements.} To determine the specific
nonconstructivities in a classical theory, and thereby to indicate
feasible directions for a constructive development, Brouwerian counterexamples
are used. The projective plane being not yet constructed here, we
give an informal example on the real metric plane, where \emph{P lies
on l} means that the distance from the point $P$ to the line $l$
is 0, while \emph{P lies outside l} means that the distance is positive.

\begin{example}
\begin{flushleft}
\label{EX. Pt and line - constructively invalid}If, on the plane
$\mathbb{R}^{2}$, we have a proof of the statement\emph{}\\
\par\end{flushleft}
\emph{\hspace{1cm}Given any point $P$ and any line $l$, either
$P$ lies on $l$, or $P$ lies outside $l$,}\\

\noindent then we have a method that will either prove the \emph{Goldbach
Conjecture, }or construct a counterexample. 
\end{example}

\begin{proof}
Using a simple finite routine, construct a sequence $\{a_{n}\}_{n\geq2}$
such that $a_{n}=0$ if $2n$ is the sum of two primes, and $a_{n}=1$
if it is not. Now apply the statement in question to the point $P=(0,\Sigma a_{n}/n^{2})$,
with the $x$-axis as the line $l$. If $P\in l$, then we have proved
the Goldbach Conjecture, while if $P\notin l$, then we have constructed
a counterexample. 
\end{proof}
For this reason, such statements are said to be \emph{constructively
invalid}. If the Goldbach question is someday settled, then other
famous problems may still be ``solved'' in this way. These examples
demonstrate that use of the \emph{Law of Excluded Middle} inhibits
mathematics from attaining its full significance. More information
concerning Brouwerian counterexamples will be found below in Section
\ref{SECTION 12. Real numbers}. 

Many other ordinary geometric statements, especially those involving
a disjunction, are also constructively invalid, admitting easily devised
Brouwerian counterexamples similar to Example \ref{EX. Pt and line - constructively invalid}.
The consequence of this Brouwerian analysis is the need for explicit
constructive details, in axioms, theorems, and proofs, which from
a classical perspective may seem superfluous. 

\emph{Constructive logic. }Following Bishop, we use no system of formal
logic. Aside from the need to avoid use of the \emph{Law of Excluded
Middle,} and to use the connective ``or\emph{''} only when warranted,
no special rules are required. The constructive logic used here is
usually called \emph{informal intuitionistic logic}; for more on this
subject, see {[}BV06, Section 1.3{]}. 

Certain concepts, such as $x=0$, for a real number $x$, are relatively
weak, compared to stronger concepts, such as $x\ne0$. The relation
$x\ne0$ requires the construction of an integer $n$ such that $1/n<|x|$;
it then follows that $x=0$ is equivalent to $\neg(x\ne0)$, while
the statement ``$\neg(x=0)$ implies $x\ne0$'' is constructively
invalid. 

In geometry, \emph{point outside a line,} $P\notin l$, is the stronger
concept, while \emph{point on a line,} $P\in l$, is the weaker.\footnote{For the stronger concept used as the single primitive notion for intuitionistic
projective geometry, see {[}D96{]}.} On the constructive real metric plane $\mathbb{R}^{2}$, the geometric
and numeric concepts are directly related; $P\notin l$ if and only
$d(P,l)>0$ {[}M07, Theorem 10.1{]}. Thus, while the statement ``If
$\neg(P\notin l)$, then $P\in l$'' will be taken as an axiom, reflecting
the constructive properties of the real numbers, the statement ``If
$\neg(P\in l)$, then $P\notin l$'' is constructively invalid. 

Further details concerning the constructive properties of the real
numbers, and constructively invalid statements, will be found below
in Section \ref{SECTION 12. Real numbers}. 

\section{{\large{}Axioms}\label{SECTION 2. Axioms}}

\noindent We adopt the usual definitions and axioms for a projective
plane, adding the several required to obtain constructive results.
The additional axioms are constructively phrased versions of elementary
facts that are immediate in classical geometry, when the \emph{Law
of Excluded Middle} is used. For a constructive study, these additional
facts must be stated explicitly in the axioms, and must be verified
whenever one constructs the finite routines for a model. 

The model $\mathbb{P}^{2}(\mathbb{R})$ in Part \ref{PART II. Analytic}
will establish the consistency of the axiom system; the question of
independence of the individual axioms is left as an open problem.
The properties of the model have served to drive the axiom choices
for the synthetic theory, taking note of Bishop's thesis, ``All mathematics
should have numerical meaning'' {[}B67, p. ix{]}.
\begin{defn}
\label{Defn. Inequality}Let $\mathscr{S}$ be a set with an equality
relation $=$. An \emph{inequality relation} $\neq$ on $\mathscr{S}$
is called an \emph{tight apartness relation}\footnote{Early work on apartness relations is due to Heyting; see {[}H66{]}.}\emph{
}if, for any $x,y,z$ in $\mathscr{S}$, the following conditions
are satisfied.

(i) $\neg(x\neq x)$.

(ii) If $x\neq y$, then $y\neq x$.

(iii) If $x\neq y$, then either $z\neq x$ or $z\neq y$.

(iv) If $\neg(x\neq y)$, then $x=y$. 
\end{defn}

Property (iii) is known as \emph{cotransitivity}, and (iv) as \emph{tightness}.
The implication ``$\neg(x=y)$ implies $x\neq y$'' is nearly always
constructively invalid, the inequality being the stronger of the two
conditions.\footnote{\noindent For a comprehensive treatment of constructive inequality
relations, see {[}BR, Section 1.2{]}.} For an example with real numbers, $x\neq0$ means that there exists
an integer $n$ such that $1/n<|x|$, while $x=0$ means only that
it is contradictory that such an integer exists.\footnote{For more details concerning the constructive properties of the real
numbers, see {[}B67, BB85, BV06{]}.} 
\begin{defn}
\noindent \label{Defn. Plane} A \emph{projective} \emph{plane $\mathbb{P}$}
consists of a family $\mathscr{P}$ of \emph{points}, and a family
$\mathscr{L}$ of \emph{lines,} satisfying the following conditions,
and axioms to be specified.

\noindent $\bullet$\emph{ Equality} relations, written =, are given
for both families $\mathscr{P}$ and $\mathscr{L}$. 

\noindent $\bullet$\emph{ Inequality} relations, written $\neq$,
with the properties of tight apartness relations, as specified in
Definition \ref{Defn. Inequality}, are given for both families $\mathscr{P}$
and $\mathscr{L}$. When $P\neq Q$, or $l\neq m$, we say that the
points \emph{$P$} and $Q$, or the lines \emph{$l$} and $m$, are
\emph{distinct}.

\noindent $\bullet$ An \emph{incidence} relation $\in$, between
the families $\mathscr{P}$ and $\mathscr{L}$, is given; when $P\in l$,
we say that the point \emph{$P$} \emph{lies on the line} $l$, and
that the line \emph{$l$} \emph{passes through} \emph{the point }$P$. 
\end{defn}

\begin{defn}
\label{Defn. outside}\emph{Outside} relation. For any point $P\in\mathscr{P}$
and any line $l\in\mathscr{L}$, we say that \emph{$P$} \emph{lies
outside the line }$l$, and that \emph{$l$} \emph{avoids the point
}$P$, written $P\notin l$, if $P\neq Q$ for all points $Q$ that
lie on $l$.
\end{defn}

A constructive definition of distinct lines is at times based on the
incidence and outside relations. Here, the relation of distinct lines
is internal, referring only to the family $\mathscr{L}$. This is
the natural approach for the model $\mathbb{P}^{2}(\mathbb{R})$ of
Part \ref{PART II. Analytic}, and is an instance where the model
influences a choice for the synthetic theory. The method here is adapted
to the situation where the families $\mathscr{P}$ and $\mathscr{L}$
are independent, as in the model, rather than the situation often
seen where lines are sets of points. With the relations of distinct
points and distinct lines established internally to the families $\mathscr{P}$
and $\mathscr{L}$, Axiom C5 will then relate the concepts to the
incidence and outside relations.\\

\noindent \textbf{Constructive Axiom Group C. }\label{Axiom-Group-C.}These
axioms form the basis for the synthetic theory. The duality of this
axiom group will be shown in Theorem \ref{Thm. Axioms dual}. 

In addition, Axiom F (Fano) will be adopted at the end of this section,
and will be shown to be self-dual in Theorem \ref{Prop. Fano self-dual};
Axiom D (Desargues) will be adopted in Section \ref{SECTION 3. Des},
where its dual (the converse) will be proved; Axiom E (Extension),
to be adopted in Section \ref{SECTION 5. Projectivities}, is self-dual;
Axiom T (the crucial component of the \emph{Fundamental Theorem}),
to be adopted in Section \ref{SECTION 6. Fund Thm}, is also self-dual,
as is Axiom P, to be adopted in Section \ref{SECTION 11. Poles and polars}
in connection with poles and polars with respect to a conic. The duality
of Definition \ref{Defn. outside}, for the \emph{outside} relation,
will be established in Theorem \ref{Thm.  outside relation dual}. 

Thus the duality of the complete set of axioms will be established.
\\

\noindent \textbf{Axiom C1.}\emph{ There exist a point $P\in\mathscr{P}$,
and a line $l\in\mathscr{L}$, such that $P\notin l$.}\\

\noindent \textbf{Axiom C2.} \emph{For any distinct points $P$ and
$Q$, there exists a unique line, denoted }$PQ$,\emph{ called the
}join, \emph{or} connecting line,\emph{ of the points, passing through
both points. }\\

\noindent \textbf{Axiom C3.}\emph{ For any distinct lines $l$ and
$m$, there exists a unique point, denoted }$l\cdot m$,\emph{ called
the }meet, \emph{or} point of intersection,\emph{ of the lines, lying
on both lines. }\\

\noindent \textbf{Axiom C4.} \emph{There exist at least three distinct
points lying on any given line. }\\

\noindent \textbf{Axiom C5.} \emph{For any lines $l$ and $m$, if
there exists a point $P\in l$ such that $P\notin m$, then $l\neq m$.}
\\

\noindent \textbf{Axiom C6.} \emph{For any point $P$ and any line
$l$, if $\neg(P\notin l)$, then $P\in l$.}\\

\noindent \textbf{Axiom C7.} \emph{If $l$ and $m$ are distinct lines,
and $P$ is a point such that $P\neq l\cdot m$, then either $P\notin l$
or $P\notin m$. }\\

\noindent \emph{Notes for Axiom Group C. }\textbf{\label{Notes - Axioms C}~ }

1. Axioms C1 thru C4 are the usual classical axioms for incidence
and extension. The remaining three axioms are statements that follow
immediately when lines are considered as sets of points, and the \emph{Law
of Excluded Middle} is used; classically, they need no explicit mention.
In this sense, no new axioms are needed for a constructive theory.

2. Axioms C2 and C3 apply only to distinct points and lines. The need
for this restriction will follow from Example \ref{EX. No Com Pt},
where it is shown that in the model $\mathbb{P}^{2}(\mathbb{R})$,
claiming the existence of a line through two arbitrary points, or
a point lying on two arbitrary lines, would be constructively invalid.

3. Axiom C6 would be immediate in a classical setting, when $P\in l$
is used in the sense of set-membership,\emph{ }where\emph{ $P\notin l$}
means \emph{$\neg(P\in l)$,} and when, applying the Law of Excluded
Middle, double negation results in an affirmative statement. 

For the constructive treatment here, the situation is quite different.
The \emph{outside} relation, \emph{$P\notin l$,} is given a strong
affirmative meaning in Definition 2.3, involving both the inequality
relation for points, and the incidence relation which connects the
two families. Just as \emph{tightness}, defined by condition (iv)
in Definition 2.1, must be assumed in Definition 2.2 for both points
and lines, the analogous condition C6, relating the two given families,
must be taken as an axiom. 

For the metric real plane, with incidence relations as noted in connection
with Example \ref{EX. Pt and line - constructively invalid}, the
condition of Axiom C6 follows from the following constructive property
of the real numbers: \emph{For any real number $\alpha$, if $\neg(\alpha\ne0)$,
then $\alpha=0$.}\footnote{This property is listed as property (i) in Section 12, where more
details are given. } For the projective model $\mathbb{P}^{2}(\mathbb{R})$, which motivates
the axiom system, Axiom C6 is verified in Corollary \ref{Cor. C6 for P2R},
using this same constructive property of the reals. 

The definitions and axioms of projective geometry may be given a wide
variety of different arrangements. For example, in {[}D96{]} the relation
$P\notin l$ is taken as a primitive notion, and the condition of
Axiom C6 becomes the definition of the incidence relation $P\in l$. 

4. Axiom C7 is a strongly worded constructive form of the statement
that distinct lines have a \emph{unique} common point. Related to
this axiom are Heyting's Axiom VI {[}H28{]}, and van Dalen's Lemma
3(f), obtained using his axiom Ax5 {[}D96{]}. Paraphrased to fit the
present context, these statements ensure that \emph{If $l$ and $m$
are distinct lines, and $P$ is a point such that $P\neq l\cdot m$
and $P\in l$, then $P\notin m$}. This is a weaker version of Axiom
C7; the stronger version will be needed here. It is an open problem
to determine whether the two versions are equivalent or constructively
distinct, or whether the weaker version would be sufficient. Generally,
a condition using the ``or'' connective is found to be constructively
stronger than other versions. 

5. An affine form of Axiom C7 is used as Axiom L1 in {[}M07{]}. 

6. Axiom C7 is the only axiom asserting a disjunction. Example \ref{EX. Pt and line - constructively invalid}
concerned the constructive invalidity of certain principles found
in classical treatments, especially those asserting a disjunction.
In Axiom C7, we have two hypotheses, each being a strong distinctness
condition. The verification of this axiom for the model $\mathbb{P}^{2}(\mathbb{R})$,
in Theorem \ref{Thm. P2R. Axioms valid.}, will require both these
strong hypotheses, other axioms, and other constructive properties
of $\mathbb{P}^{2}(\mathbb{R})$. 

7. Axiom C7 may rightly claim a preëminent standing in the axiom system;
it will be indispensable for nearly all the constructive proofs. 
\begin{prop}
\label{Prop. lines to two points}Let $P,Q,R$ be distinct points.
Then $P\notin QR$ if and only if $PQ\neq PR$. 
\end{prop}

\begin{proof}
First let $P\notin QR$. From Axiom C5 we have $PR\neq QR$. Since
$Q\neq R=PR\cdot QR$, it follows from Axiom C7 that $Q\notin PR$;
thus $PQ\neq PR$. Conversely, if $PQ\neq PR$, then from $P\neq Q=PQ\cdot QR$
it follows that $P\notin QR$. 
\end{proof}

\begin{prop}
\label{Prop. Converse to Axiom C5}If the lines $l$ and $m$ are
distinct, then there exists a point $P\in l$ such that $P\notin m$. 
\end{prop}

\begin{proof}
Set $Q=l\cdot m$, using Axiom C3, and select a point $P\in l$ such
that $P\neq Q$, using Axiom C4. It follows from Axiom C7 that either
$P\notin l$ or $P\notin m$. The first case is ruled out by Axiom
C6; thus $P\notin m$. 
\end{proof}
\begin{defn}
\label{Defn. Collinear etc etc}~~

\noindent \emph{$\bullet$} A set $\mathscr{S}$ of points is \emph{collinear}
if $P\in QR$ whenever $P,Q,R\in\mathscr{S}$ with $Q\neq R$. 

\noindent \emph{$\bullet$} A set $\mathscr{S}$ of points is \emph{noncollinear}
if there exist distinct points $P,Q,R\in\mathscr{S}$ such that $P\notin QR$.

\noindent $\bullet$ A set $\mathscr{T}$ of lines is \emph{concurrent}
if ~$l\cdot m\in n$ whenever $l,m,n\in\mathscr{T}$ with $l\neq m$. 

\noindent \emph{$\bullet$} A set $\mathscr{T}$ of lines is \emph{nonconcurrent}
if there exist distinct lines $l,m,n\in\mathscr{T}$ such that $l\cdot m\notin n$.

\noindent $\bullet$ The \emph{range} \emph{of points on a line} \emph{$l$}
is the set $\overline{l}=\{P\in\mathscr{P}:P\in l\}$. 

\noindent $\bullet$ The \emph{pencil} \emph{of lines} \emph{through
a point} $Q$ is the set $Q^{*}=\{m\in\mathscr{L}:Q\in m\}$.\emph{ }
\end{defn}

\begin{example}
\noindent \label{Ex. collinear}A stronger, classically equivalent,
alternative definition for \emph{collinear set} is the condition \emph{There
exists a line that passes through each point of the set}. The equivalence
of the two conditions is constructively invalid for the model $\mathbb{P}^{2}(\mathbb{R})$.
The Brouwerian counterexample given in Example \ref{EX. No Com Pt}
will apply; we give a simplified version here, in brief form. Consider
two points on the plane $\mathbb{R}^{2}$, the origin, and a point
close to or at the origin. The set formed by these has at most two
points, and is collinear according to our definition, yet it is not
possible, constructively, to predict what line might contain both
points.\footnote{For more on the constructive eccentricities of such sets, see {\small{}{[}M13b,
}Example 2.5{\small{}{]}}. }

Another distinction between the alternative definitions concerns the
statement \emph{If  $\neg(\mathscr{S}$ is noncollinear), then $\mathscr{S}$
is collinear.} This statement follows easily from our definition,
while under the alternative definition it is seen to be constructively
invalid, using the example above. 
\end{example}

\begin{prop}
\noindent If $\mathscr{S}$ is a noncollinear set of points, then
for any line $l$ in the plane, there exists a point in $\mathscr{S}$
that lies outside $l$. 
\end{prop}

\begin{proof}
Choose distinct points $P,Q,R\in\mathscr{S}$ as in Definition \ref{Defn. Collinear etc etc},
with $P\notin QR$. It follows from Proposition \ref{Prop. lines to two points}
that $PQ\neq PR$. By cotransitivity for lines, either $l\neq PQ$
or $l\neq PR$. It suffices to consider the first case; set $Y=l\cdot PQ$.
Now, either $Y\neq P$ or $Y\neq Q$. In the first subcase, we have
$P\neq l\cdot PQ$, so it follows from Axiom C7 that $P\notin l$.
Similarly, in the second subcase we find that $Q\notin l$. 
\end{proof}
\begin{prop}
\noindent \label{Prop. 3 noncoll points}If a set $\mathscr{S}$ of
three distinct points is noncollinear, then $P\notin QR$, where $P,Q,R$
are the points of \textup{$\mathscr{S}$} taken in any order. 
\end{prop}

\begin{proof}
\noindent Given that $P\notin QR$, we have $PR\neq QR$ and $Q\neq R=PR\cdot QR$;
thus by Axiom C7 it follows that $Q\notin PR$. By symmetry, we also
have $R\notin PQ$. 
\end{proof}
Given any statement, the \emph{dual statement} is obtained by interchanging
the words ``point'' and ``line''.
\begin{thm}
\label{Thm. Axioms dual}\emph{.} The definition of the projective
plane $\mathbb{P}$ is self-dual, and the dual of each axiom in Axiom
Group C holds on $\mathbb{P}$. 
\end{thm}

\begin{proof}
Definition \ref{Defn. Plane}, and Axioms C1, C2/C3, C6 are clearly
self-dual. 

For the dual of Axiom C4, select a point $Q$ and a line $m$, with
$Q\notin m$, using Axiom C1. Using Axiom C4 select three distinct
points $R_{1},R_{2},R_{3}$ on $m$. Using Definition \ref{Defn. outside},
we have $Q\neq R_{i}$ for each $i$; set $l_{i}=QR_{i}$. Since $Q\notin m$,
it follows that $l_{i}\neq m$ for each $i$. Since $R_{1}\neq R_{2}=m\cdot l_{2}$,
it follows from Axiom C7 that $R_{1}\notin l_{2}$, so $l_{1}\neq l_{2}$.
By symmetry, the three lines $l_{i}$ are distinct. This is the desired
result for the selected point $Q$, based on the existence of the
line $m$ that avoids $Q$. Now, given an arbitrary point $P$, using
cotransitivity we may assume that $P\neq R_{1}$. Since $P\neq R_{1}=m\cdot l_{1}$,
it follows that either $P\notin m$ or $P\notin l_{1}$. In either
case, using the same method as for $Q$ and $m$, we may construct
three distinct lines through $P$. 

The dual of Axiom C5 states that \emph{``Points $P$ and $Q$ are
distinct, $P\neq Q$, if and only if there exists a line $l$ such
that $P\in l$ and $Q\notin l$, or vice-versa''}. First, consider
points $P$ and $Q$ with $P\neq Q$, and use Axiom C4(dual) to construct
distinct lines $l$ and $m$ through $P$. Now we have $Q\neq P=l\cdot m$,
so by Axiom C7 we may assume that $Q\notin l$. Thus we have a line
$l$ through $P$ that avoids $Q$. Conversely, if for some line $l$
we have $Q\notin l$ and $P\in l$, then $Q\neq P$ by Definition
\ref{Defn. outside}.

The dual of Axiom C7 states that ``\emph{If $Q$ and $R$ are distinct
points, and $n$ is a line such that $n\neq QR$, then either $Q\notin n$
or $R\notin n$''. }To prove this, set $S=QR\cdot n$, and use cotransitivity
to obtain either $S\neq Q$ or $S\neq R$. In the first case, since
$Q\neq S=QR\cdot n$, it follows from Axiom C7 that $Q\notin n$.
Similarly, in the second case we find that $R\notin n$. Thus Axiom
C7 is self-dual. 
\end{proof}
\begin{thm}
\label{Thm.  outside relation dual}Let $P$ be any point, and $l$
any line. Then $P\notin l$ if and only if \emph{$l\neq m$ for all
lines $m$ that pass through $P$. }
\end{thm}

\begin{proof}
\noindent Let $P\notin l$ and let $m$ be any line through $P$.
Using Axiom C4(dual) and cotransitivity for lines, select a line $n$
passing through $P$ and distinct from $m$, and select any point
$Q\in l$. Then $Q\neq P=m\cdot n$, so by Axiom C7 it follows that
either $Q\notin m$, or $Q\notin n$. In the first case, $l\neq m$.
In the second case, $l\neq n$; set $R=l\cdot n$. Since $R\in l$,
we have $P\neq R$; thus $R\neq P=m\cdot n$, so $R\notin m$, and
again $l\neq m$. Thus the dual condition is satisfied.

Now let $P$ and $l$ satisfy the dual condition, and let $Q$ be
any point on $l$. Select a point $R$ on $l$, distinct from $Q$.
Either $P\neq Q$ or $P\neq R$. In the second case, set $m=PR$;
by hypothesis, $l\neq m$. Since $Q\neq R=l\cdot m$, it follows that
$Q\notin m$. From Axiom C5(dual), we have $P\neq Q$. Hence $P\notin l$. 
\end{proof}
\begin{cor}
\label{Cor-primary-relation-self-dual}The primary relation adopted
in Definition \ref{Defn. outside}, \emph{point outside a line}, is
self-dual in the context of Axiom Group C. 
\end{cor}

\noindent From this corollary, and Theorem \ref{Thm. Axioms dual},
we obtain the \emph{duality principle:}
\begin{thm}
\label{Thm-Duality}On the projective plane $\mathbb{P}$, the dual
of any result is immediately valid, with no further proof required. 
\end{thm}

\begin{defn}
\label{Defn. Quadrangle}A \emph{quadrangle} is an ordered set $PQRS$
of four distinct points, the \emph{vertices, }such that each subset
of three points is noncollinear. The \emph{sides} are the six lines
joining the vertices. The three \emph{diagonal points} are $D_{1}=PQ\cdot RS$,
$D_{2}=PR\cdot QS$, and $D_{3}=PS\cdot QR$. A \emph{quadrilateral},
with four sides and six vertices, is the dual configuration. 
\end{defn}

\noindent \emph{Note for Definition \ref{Defn. Quadrangle}.} Since
$R\notin PQ$, we have $PQ\neq RS$; by symmetry, all six sides are
distinct, and the definition of the diagonal points is valid. By cotransitivity,
either $D_{2}\neq P$ or $D_{2}\neq R$. It suffices to consider the
first case; thus $D_{2}\neq PQ\cdot PR$, so by Axiom C7 we have $D_{2}\notin PQ$,
and $D_{2}\neq D_{1}$. By symmetry, all three diagonal points are
distinct.\\

It will be convenient to exclude certain finite planes, such as the
seven-point ``Fano plane'' {[}Fan92{]}, an illustration of which
may be found at {[}VY10, p. 45{]} or {[}Wei07, p. 1294{]}. Thus we
adopt the following:\\

\noindent \textbf{Axiom F.\label{Axiom F. Fano.}} Fano's Axiom. \emph{The
diagonal points of any quadrangle are noncollinear.}
\begin{prop}
\label{Prop. Fano self-dual}Axiom F is self-dual; the diagonal lines
of any quadrilateral are nonconcurrent. 
\end{prop}

\begin{proof}
Given a quadrilateral $pqrs$, denote four of the six vertices as
$P=p\cdot q$, $Q=q\cdot r$, $R=r\cdot s$, $S=s\cdot p$. The diagonal
lines of $pqrs$ are then $d_{1}=(p\cdot q)(r\cdot s)=PR$, $d_{2}=(p\cdot r)(q\cdot s)$,
and $d_{3}=(p\cdot s)(q\cdot r)=QS$. 

To show that $PQRS$ is a quadrangle, we note that since the lines
$p,q,r$ are nonconcurrent, we have $P=p\cdot q\notin r=QR$, so the
points $P,Q,R$ are noncollinear, and similarly for the other three
triads. 

The diagonal points of the quadrangle $PQRS$ are $D_{1}=PQ\cdot RS=q\cdot s$,
$D_{2}=PR\cdot QS=d_{1}\cdot d_{3}$, and $D_{3}=PS\cdot QR=p\cdot r$;
it follows that $d_{2}=D_{1}D_{3}$. By Axiom F, we have $D_{2}\notin D_{1}D_{3}$;
thus $d_{1}\cdot d_{3}\notin d_{2}$, and the diagonal lines $d_{1},d_{2},d_{3}$
of the quadrilateral are nonconcurrent. 
\end{proof}

\section{{\large{}Desargues's Theorem} \label{SECTION 3. Des}}

\noindent We adopt Desargues's Theorem as Axiom D, and then use it
to prove the converse, which is its dual.\footnote{Desargues's Theorem and the converse are both derived in {[}H28, §§5-6{]},
using axioms for projective space; here we use only axioms for a plane. }
\begin{defn}
~\label{Defn. Triangle}

\noindent $\bullet$ A \emph{triangle} is an ordered triad \emph{$PQR$}
of distinct, noncollinear points. The three points are the \emph{vertices};
the lines $PQ$,$PR$,$QR$ are the \emph{sides}. 

\noindent $\bullet$ Triangles \emph{$PQR$} and \emph{$P'Q'R'$ }are
\emph{distinct} if corresponding vertices are distinct, and corresponding
sides are distinct.

\noindent $\bullet$ Distinct triangles are said to be \emph{perspective
from the center $O$} if the three lines joining corresponding vertices
are concurrent at $O$, and $O$ lies outside each of the six sides. 

\noindent $\bullet$ Distinct triangles are said to be \emph{perspective
from the axis $l$} if the three points of intersection of corresponding
sides are collinear on $l$, and $l$ avoids each of the six vertices.
\end{defn}

\noindent \begin{flushleft}
\textbf{Axiom D. }\label{Ax-D. Des}\emph{If two triangles are perspective
from a center, then they are also perspective from an axis. }
\par\end{flushleft}
\begin{thm}
\label{Thm. DesConv}If two triangles are perspective from an axis,
then they are also perspective from a center.
\end{thm}

\begin{proof}
We are given distinct triangles $PQR$, $P'Q'R'$, with points $A=QR\cdot Q'R'$,
$B=PR\cdot P'R'$, $C=PQ\cdot P'Q'$ collinear on a line $l$, with
$V\notin l$ for all six vertices $V$. 

Since $Q\notin l$, we have $A\neq Q=PQ\cdot QR$; it follows from
Axiom C7 that $A\notin PQ$, and thus $A\neq C$. By symmetry, all
three points $A,B,C$ are distinct, the points $A,Q,Q'$ are distinct,
and the points $B,P,P'$ are distinct. Since $Q\neq A=QR\cdot Q'R'$,
we have $Q\notin Q'R'=AQ'$; thus the points $A,Q,Q'$ are noncollinear,
and similarly for $B,P,P'$. Since $P\notin l=AB$, it follows that
$AB\neq BP$. Since $A\neq B=AB\cdot BP$, we have $A\notin BP$,
so $AQ\neq BP$; similarly, $AQ'\neq BP'$. Since $Q'\neq C=PQ\cdot P'Q'$,
it follows that $Q'\notin PQ$; thus $QQ'\neq PQ$, and similarly
$PP'\neq PQ$. Since $P\neq Q=PQ\cdot QQ'$, we have $P\notin QQ'$;
thus $PP'\neq QQ'$. Set $O=PP'\cdot QQ'$. 

The above shows that the auxiliary triangles $AQQ'$, $BPP'$ are
distinct. The lines $AB$, $PQ$, $P'Q'$, joining corresponding vertices,
are concurrent at $C$. Since $C\neq A=l\cdot AQ$, it follows that
$C\notin AQ$; similarly, $C\notin AQ'$. Since $P\notin QQ'$, it
follows that $CQ=CP\neq QQ'$, and from $C\neq Q=CQ\cdot QQ'$ we
have $C\notin QQ'$. Thus $C$ lies outside each side of triangle
$AQQ'$, and similarly for triangle $BPP'$. 

Thus the auxiliary triangles $AQQ'$, $BPP'$ are perspective from
the center $C$; it follows from Axiom D that these triangles are
perspective from the axis $(AQ\cdot BP)(AQ'\cdot BP')=RR'$. Thus
$O\in RR'$, and the axis $RR'$ avoids all six vertices of the auxiliary
triangles. This shows that the lines $PP'$, $QQ'$, $RR'$ are concurrent
at $O$. Since $P\notin RR'$, we have $P\neq O$. Since $O\neq P=PP'\cdot PQ$,
it follows that $O\notin PQ$. Similarly, $O$ lies outside each side
of the triangles $PQR$, $P'Q'R'$. Hence the original triangles are
perspective from the center $O$. 
\end{proof}

\section{{\large{}Harmonic conjugates}\label{SECTION 4. HarConj} }

\noindent Harmonic conjugates are often defined using quadrangles
or triangles. We must use a less problematic definition; it must apply
to every point on the line, including the base points, and any point
for which it is not known, constructively, whether or not it coincides
with a base point. The simplicity of the definition will facilitate
the verification that the joins and intersections used are constructively
meaningful, and that the result is independent of the selection of
construction elements. 
\begin{defn}
\label{Defn. HarConj}Let\emph{ $A$} and \emph{$B$} be distinct
points. For any point \emph{$C$} on the line $AB$, select a line
\emph{$l$} through $C$, distinct from $AB$, and select a point
\emph{$R$} lying outside each of the lines \emph{$AB$} and $l$.
Set $P=BR\cdot l$, $Q=AR\cdot l$, and \emph{$S=AP\cdot BQ$}. Pending
verifications in Proposition \ref{Prop. HC-verify} and Theorem \ref{Thm. HC-indep},
the point \emph{$D=AB\cdot RS$} will be called the \emph{harmonic
conjugate} \emph{of} $C$ \emph{with respect to the points} $A,B$;
we write $D=h(A,B;C)$. The points $A,B,C,D$ are said to form a \emph{harmonic
set}, written $h(A,B;C,D)$. 
\end{defn}

\begin{lem}
\label{Lm. HarCon.details}In Definition \ref{Defn. HarConj}: 

\emph{(a)} $P\neq A$, $Q\neq B$, $P\neq Q$.

\emph{(b)} $P\notin AR$, $Q\notin BR$, $A\notin BR$, $B\notin AR$.

\emph{(c)} $AR\neq BR,$ $AP\neq AR$, $AP\neq BR$, $BQ\neq BR$,
$BQ\neq AR$.
\end{lem}

\begin{proof}
Since $A\neq B=AB\cdot BR$, it follows from Axiom C7 that $A\notin BR$;
thus $A\neq P$, $AP\neq BR$, and $AR\neq BR$. By symmetry, $B\notin AR$,
$B\neq Q$, and $BQ\neq AR$. Since $P\neq R=AR\cdot BR$, we have
$P\notin AR$, so $P\neq Q$, and $AP\neq AR$. Similarly, $Q\notin BR$
and $BQ\neq BR$. 
\end{proof}

\begin{prop}
\noindent \label{Prop. HC-verify}The construction of a harmonic conjugate,
in Definition \ref{Defn. HarConj}\emph{, }involves valid joins and
intersections.
\end{prop}

\begin{proof}
Using Lemma \ref{Lm. HarCon.details}, we need only show that $AP\neq BQ$,
and that $R\neq S$. By cotransitivity, we may assume that $C\neq B$.
Since $C\neq B=AB\cdot BR$, it follows from Axiom C7 that $C\notin BR$,
so $C\neq P$. Since $P\neq C=AB\cdot l$, it follows that $P\notin AB$;
thus $AB\neq AP$. Since $B\neq A=AB\cdot AP$, we have $B\notin AP$;
thus $AP\neq BQ$. This shows that the definition of $S$ is valid.
From $R\neq B=BQ\cdot BR$, it follows that $R\notin BQ$; hence $R\neq S$.
This shows that the definition of $D$ is valid. 
\end{proof}
The next result is one of the four lemmas required for the proof of
Theorem \ref{Thm. HC-indep}, which will validate the harmonic conjugate
construction, showing that it is independent of the choice of construction
elements. This lemma involves the special situation in which the original
point is one of the base points. 
\begin{lem}
\label{Lm. HC. C =00003D A}Let $A\neq B$. Then $h(A,B;A)=A$ and
$h(A,B;B)=B$, for any selection of construction elements $l,R$ in
Definition \ref{Defn. HarConj}. 
\end{lem}

\begin{proof}
When $C=A$, then $Q=AR\cdot l=CR\cdot l=C=A$, so $S=AP\cdot BQ=AP\cdot BA=A$,
and thus $D=AB\cdot RS=AB\cdot RA=A$. Similarly when $C=B$.
\end{proof}
\begin{lem}
\label{Lm. C not A. C not B.}In Definition \ref{Defn. HarConj}: 

\emph{(a)} If $C\neq A$, then $Q\notin AB$, $Q\neq S$, $S\neq A$,
and $D\neq A$.

\emph{(b)} If $C\neq B$, then $P\notin AB$, $P\neq S$, $S\neq B,$
and $D\neq B$.
\end{lem}

\begin{proof}
It will suffice to consider (a), since (b) will follow by symmetry.
Since $C\neq A=AB\cdot AR$, it follows from Axiom C7 that $C\notin AR$;
thus $C\neq Q$. Since $Q\neq C=AB\cdot l$, we have $Q\notin AB$,
so $Q\neq A$. Since $A\neq Q=BQ\cdot AR$, it follows that $A\notin BQ$;
thus $A\neq S$. Since $S\neq A=AP\cdot AR$, we have $S\notin AR$,
so $S\neq Q$, and $AR\neq RS$. Since $A\neq R=AR\cdot RS$, it follows
that $A\notin RS$, and thus $A\neq D$. 
\end{proof}
The next lemma shows that for a point distinct from both base points,
the traditional quadrangle will appear; for a complete statement regarding
this configuration, see Corollary \ref{Cor. HC-Quad.}.
\begin{lem}
\label{Lm. C not ei base point}In Definition \ref{Defn. HarConj}
for the construction of a harmonic conjugate, let $C\neq A$ and $C\neq B$.
Then the four points $P,Q,R,S$ are distinct and lie outside the line
$AB$, and each subset of three points is noncollinear. Furthermore,
$h(A,B;C)\neq C$.
\end{lem}

\begin{proof}
Using Lemma \ref{Lm. C not A. C not B.}, we see that $P,Q,R$ lie
outside $AB$. By the same lemma, we also have $S\neq B=AB\cdot BQ$,
so it follows from Axiom C7 that $S\notin AB$. Thus the four points
lie outside $AB$.

In Definition \ref{Defn. HarConj}, we have $P\neq R$, $Q\neq R$.
From Lemma \ref{Lm. HarCon.details}, we have $P\neq Q$. From Proposition
\ref{Prop. HC-verify}, we have $R\neq S$. From Lemma \ref{Lm. C not A. C not B.},
we have $P\neq S$, $Q\neq S$. Thus the four points are distinct. 

For the triads, we use the results of Lemmas \ref{Lm. HarCon.details}
and \ref{Lm. C not A. C not B.}. Since $S\neq A=AP\cdot AR$, it
follows that $S\notin AR$, so $AR\neq RS$. Since $Q\neq R=AR\cdot RS$,
it follows that $Q\notin RS$; thus the points $Q,R,S$ are noncollinear,
and similarly for $P,R,S$. Since $P\neq B=BQ\cdot BR$, we have $P\notin BQ=QS$;
thus the points $P,Q,S$ are noncollinear. Since $P\notin AR=QR$,
the points $P,Q,R$ are noncollinear. Thus each triad is noncollinear. 

Since two diagonal points of the quadrangle $PQRS$ are $A=PS\cdot QR$
and $B=PR\cdot QS$, Axiom F, Fano's Axiom, asserts that the third
diagonal point $T$ lies outside $AB$. Thus $C\neq T=PQ\cdot RS$;
it follows that $C\notin RS$, and hence $C\neq D$.
\end{proof}
The next theorem will validate the harmonic conjugate construction.\footnote{A harmonic conjugate construction based on quadrangles is validated
in {[}H28, §7{]}, using axioms for projective space; the construction
applies only to points distinct from the base points. Here we use
only axioms for a plane, and consider all points on the base line. }
\begin{thm}
\label{Thm. HC-indep}The construction of the harmonic conjugate of
an arbitrary point C on a line AB, with respect to the points A,B,
using Definition \ref{Defn. HarConj}, results in a point D that is
independent of the selections of the line l and the point R used in
the construction.
\end{thm}

\begin{proof}
Let $l',R'$ be alternative selections, and let $D'$ be the point
resulting when the alternatives are used in the construction. 

(1) Suppose that $D'\neq D$, and suppose further that $C\neq A$,
$C\neq B$, $l'\neq l$, and $R'\neq R$. We will contradict these
five assumptions sequentially, in reverse. This process ends in a
negation of the first assumption; thus the required conclusion, $D'=D$,
will follow from the tightness property of the inequality relation
for points, specified in Definition \ref{Defn. Inequality}(iv). 

(2) Using only the first four assumptions in (1), we note here a few
basic facts. By Lemma \ref{Lm. C not ei base point}, the points $P,Q,R,S$
are distinct and lie outside the line $AB$, and each subset of three
points is noncollinear. Similarly for the points $P',Q',R',S'$. 

Since $D\neq D'=AB\cdot R'S'$, it follows from Axiom C7 that $D\notin R'S'$;
thus $RS\neq R'S'$. From $P\neq C=l\cdot l'$ it follows that $P\notin l'$.
Thus $P\neq P'$, and similarly, $Q\neq Q'$. Since $P\neq Q$, we
have $PQ=l$; similarly, $P'Q'=l'$. Thus $PQ\neq P'Q'$. From $Q'\neq C=l\cdot l'$
it follows that $Q'\notin l=PQ$, and thus $QQ'\neq PQ$. Since $P\neq Q=QQ'\cdot PQ$,
we have $P\notin QQ'$, so $PP'\neq QQ'$. Thus we may define $O=PP'\cdot QQ'$. 

(3) Since $R'\neq R=AR\cdot BR$, it follows that either $R'\notin AR$
or $R'\notin BR$. By symmetry, it suffices to consider the second
case. Since $PR=BR$, it follows that $PR\ne P'R'$.

(4) Suppose further, in addition to the assumptions in (1), that $R'\notin AR$,
$PS\neq P'S'$, $QS\neq Q'S'$, $O\notin RS$, $O\notin R'S'$, and
$S\neq S'$. 

(5)\emph{ }Consider the triangles $PQR$, $P'Q'R'$. From (2) and
the fifth assumption in (1), the corresponding vertices of the triangles
are distinct, and $PQ\neq P'Q'$. Since $R'\notin AR$ by (4), it
follows that $QR=AR\neq AR'=Q'R'$. Also, $PR\ne P'R'$ by (3). Thus
the triangles are distinct. Since $QR\cdot Q'R'=A$, $PR\cdot P'R'=B$,
$PQ\cdot P'Q'=C$, and the six vertices lie outside the line $AB$,
these triangles are perspective from the axis $AB$. By the converse
to Desargues's Theorem, established above as Theorem \ref{Thm. DesConv},
the triangles are perspective from the center $PP'\cdot QQ'=O$. 

Thus $O\in RR'$, and $O$ lies outside each of the six sides of the
triangles $PQR$, $P'Q'R'$. 

(6) The triangles $PQS$, $P'Q'S'$ are distinct, using the assumptions
in (4), with $PS\cdot P'S'=A$, $QS\cdot Q'S'=B$, $PQ\cdot P'Q'=C$,
and vertices outside the line $AB$. Thus these triangles are perspective
from the axis $AB$; it follows that they are perspective from the
center $PP'\cdot QQ'=O$. 

Thus $O\in SS'$, and $O$ lies outside each of the six sides of the
triangles $PQS$, $P'Q'S'$.

(7) Now consider the triangles $PRS$, $P'R'S'$. By (1, 2, 4), the
corresponding vertices are distinct. From (2) we have $RS\neq R'S'$,
from (3) we have $PR\ne P'R'$, and in (4) we assumed that $PS\neq P'S'$.
Thus the triangles are distinct. From (5) we have $O\notin PR$, from
(6) we have $O\notin PS$, and from (4) we have $O\notin RS$. Similarly,
$O$ lies outside each side of the second triangle. By (5) and (6),
the triangles are perspective from the center $O$; thus by Desargues's
Theorem, adopted above as Axiom D, they are also perspective from
the axis $(PS\cdot P'S')(PR\cdot P'R')=AB$. Thus the point $E:=RS\cdot R'S'$
lies on $AB$. It follows that $D=E$, and also $D'=E$. This contradiction
of the first assumption in (1) negates the last assumption in (4). 

Thus $S=S'$.

(8) It follows that $BS=BS'$; i.e., $QS=Q'S'$. This contradiction
of the third assumption in (4) negates the fifth; thus $O\in R'S'$. 

Similarly, we obtain $PS=P'S'$, a contradiction of the second assumption
in (4), negating the fourth; thus $O\in RS$. The condition $QS=Q'S'$
negates the third assumption in (4); hence $QS=Q'S'$, and thus $S\in QQ'$.
The second assumption in (4) is negated similarly; thus $PS=P'S'$,
and $S\in PP'$. 

(9) From (8), it is easily seen (in two ways) that $S=O$. In step
(5), the arguments depend only on (1) and the first assumption in
(4), not yet negated; thus we may use here the conclusion that $O\in RR'$.
It follows that $RS=R'S'$, a contradiction of (2), negating the first
assumption in (4). 

Thus $R'\in AR$, and $AR=AR'$. 

(10) From (8, 9), we have $QS\cdot AR=Q'S'\cdot AR'$, i.e., $Q=Q'$,
contradicting (2), and negating the last assumption in (1). Thus $R'=R$.
Combining this with (7), we have $RS=R'S'$. This contradiction of
(2) negates the fourth assumption in (1). Thus $l'=l$. 

(11) By (10), it is evident that $D'=D$, contradicting the first
assumption in (1), and negating the third; thus $C=B$. By Lemma \ref{Lm. HC. C =00003D A},
it follows that $D=B$, and also $D'=B$; this contradicts the first
assumption in (1), and negates the second. Thus $C=A$; using the
same lemma, this again results in a contradiction, negating the first
assumption in (1). 

Hence $D'=D$, and this validates the harmonic conjugate construction. 
\end{proof}
\begin{cor}
\label{Cor. HC-Quad.}Let $A,B,C,D$ be collinear points, with $C$
distinct from both points $A$ and $B$. Then $D=h(A,B;C$) if and
only if there exists a quadrangle $PQRS$, with vertices outside the
line $AB$, such that $A=PS\cdot QR$, $B=QS\cdot PR$, $C\in PQ$,
and $D\in RS$. 
\end{cor}

\begin{lem}
\label{Lm. hc of D is C}Let $A\neq B$, and let $C$ and $D$ be
any points on the line $AB$.

\emph{(a)} If $h(A,B;C,D)$, then $h(B,A;C,D)$.

\emph{(b)} If $D=h(A,B;C)$, then $C=h(A,B;D)$\textup{.} 
\end{lem}

\begin{proof}
(a) This follows from the symmetry of the construction in Definition
\ref{Defn. HarConj}.

(b) By cotransitivity and (a), it suffices to consider the case in
which $C\neq A$. Using the notation of Definition \ref{Defn. HarConj}
for the construction of $D$, from Lemma \ref{Lm. C not A. C not B.}
we have $Q\notin AB$, $Q\neq S$, $S\neq A$, and $D\neq A$. Define
$l^{d}=RS$ and $R^{d}=Q$. Then $D\in l^{d}$, $l^{d}\neq AB$, and
$R^{d}\notin AB$. Since $D\neq A=AB\cdot AR$, it follows from Axiom
C7 that $D\notin AR$, and thus $AR\neq RS$. Since $Q\neq R=AR\cdot RS$,
we have $Q\notin RS$; i.e., $R^{d}\notin l^{d}$. Thus the elements
$l^{d}$, $R^{d}$ may be used to construct $h(A,B;D)$. 

Now $P^{d}=BR^{d}\cdot l^{d}=BQ\cdot RS=S$, and $Q^{d}=AR^{d}\cdot l^{d}=AQ\cdot RS=AR\cdot RS=R$.
Thus $S^{d}=AP^{d}\cdot BQ^{d}=AS\cdot BR=AP\cdot BR=P$. It follows
that $R^{d}S^{d}=QP=l$ and $h(A,B;D)=AB\cdot R^{d}S^{d}=AB\cdot l=C$.
\end{proof}
\begin{thm}
\label{Thm. harcon.bij}Let $A$ and $B$ be distinct points in a
range $\overline{l}$, and let $\upsilon$ be the mapping of harmonic
conjugacy with respect to the base points $A,B$; i.e., set $X^{\upsilon}=h(A,B;X)$,
for all points $X$ in $\overline{l}$. Then $\upsilon$ is a bijection
of the range $\overline{l}$ onto itself, of order 2.
\end{thm}

\begin{proof}
Lemma \ref{Lm. hc of D is C}(b) shows that $\upsilon$ is onto $\overline{l}$,
and of order 2. Now let $C_{1},C_{2}\in AB$, with $C_{1}\neq C_{2}$.
To show that the harmonic conjugates $D_{1}=h(A,B;C_{1})$ and $D_{2}=h(A,B;C_{2})$
are distinct, we note first that by cotransitivity either $C_{1}\neq A$
or $C_{1}\neq B$, and similarly for $C_{2}$. By symmetry, only two
of the four cases need be considered. 

\emph{Case 1;} $C_{1}\neq A$ \emph{and} $C_{2}\neq A$\emph{.} Select
a point $R$ with $R\notin AB$, and select a point $Q\in AR$ with
$Q\neq A$ and $Q\neq R$. Since $Q\neq A=AB\cdot AR$, it follows
from Axiom C7 that $Q\notin AB$; thus $Q\neq C_{1}$ and $Q\neq C_{2}$.
Set $l_{1}=C_{1}Q$ and $l_{2}=C_{2}Q$; thus $l_{1}\neq AB$ and
$l_{2}\neq AB$. Since $A\neq C_{1}=AB\cdot l_{1}$, it follows that
$A\notin l_{1}$, so $AR\neq l_{1}$; similarly, $AR\neq l_{2}$.
Since $C_{1}\neq C_{2}=AB\cdot l_{2}$, we have $C_{1}\notin l_{2}$,
so $l_{1}\neq l_{2}$. Since $R\neq Q=AR\cdot l_{1}$, it follows
that $R\notin l_{1}$; similarly, $R\notin l_{2}$. Thus $l_{1},R$
and $l_{2},R$ may be used in Definition \ref{Defn. HarConj} to construct
the harmonic conjugates $D_{1}$ and $D_{2}$. 

Clearly, $Q_{1}=Q_{2}=Q$. Since $P_{1}\neq Q=l_{1}\cdot l_{2}$,
it follows that $P_{1}\notin l_{2}$, and thus $P_{1}\neq P_{2}$.
Since $P_{2}\neq P_{1}=AP_{1}\cdot BR$, we have $P_{2}\notin AP_{1}$,
so $AP_{1}\neq AP_{2}$. Since $S_{1}\neq A=AP_{1}\cdot AP_{2}$,
it follows that $S_{1}\notin AP_{2}$, and thus $S_{1}\neq S_{2}$.
Since $R\neq Q=AR\cdot BQ$, we have $R\notin BQ$, so $BQ\neq RS_{2}$.
Since $S_{1}\neq S_{2}=BQ\cdot RS_{2}$, it follows that $S_{1}\notin RS_{2}$,
and thus $RS_{1}\neq RS_{2}$. Since $D_{1}\neq R=RS_{1}\cdot RS_{2}$,
we have $D_{1}\notin RS_{2}$, and hence, finally, $D_{1}\neq D_{2}$.

\emph{Case 2;} $C_{1}\neq A$ \emph{and} $C_{2}\neq B$\emph{.} By
Lemma \ref{Lm. C not A. C not B.} it follows that $D_{1}\neq A$.
From cotransitivity it follows that either $D_{2}\neq D_{1}$ or $D_{2}\neq A$,
so we may assume that $D_{2}\neq A$. From Lemmas \ref{Lm. C not A. C not B.}
and \ref{Lm. hc of D is C}(b), we have $C_{2}\neq A$, and now Case
1 applies. 
\end{proof}

\section{{\large{}Projectivities} \label{SECTION 5. Projectivities} }

\noindent We use the Poncelet {[}Pon22{]} definition of \emph{projectivity}.
Theorem \ref{Thm, projty-preserves-harcons} will show that every
Poncelet projectivity is a von Staudt {[}Sta47{]} projectivity. 
\begin{defn}
\noindent \label{Defn. Projection}~

\noindent $\bullet$ The \emph{projection} $\rho:\overline{l}\rightarrow\overline{m}$,
of a range $\overline{l}$ of points onto a range $\overline{m}$,
from the\emph{ center} $T$, where the point $T$ lies outside both
lines $l$ and $m$, is the bijection defined by $X^{\rho}=TX\cdot m$,
for all points $X$ in the range $\overline{l}$. We write $\rho=\rho(T;l,m)$. 

\noindent $\bullet$ The \emph{projection} $\rho:P^{*}\rightarrow Q^{*},$
of a pencil $P^{*}$ of lines onto a pencil $Q^{*}$, by the \emph{axis}
$n$, where the line $n$ avoids both points $P$ and $Q$, is the
bijection defined by $l^{\rho}=(n\cdot l)Q$, for all lines $l$ in
the pencil $P^{*}$. We write $\rho=\rho(n;P,Q)$. 

\noindent $\bullet$ The \emph{section} of a pencil $P^{*}$, by a
line $m$ that avoids the point $P$, is the bijection $\rho:P^{*}\rightarrow\overline{m}$
defined by $l^{\rho}=l\cdot m$, for all lines $l$ in the pencil
$P^{*}$. The dual and inverse $\rho^{-1}:\overline{m}\rightarrow P^{*}$,
defined by $X\rightarrow PX$, for all points $X$ in the range $\overline{m}$,
is also called a \emph{section}. 

\noindent $\bullet$ Any projection or section is said to be a \emph{perspectivity}.
\end{defn}

\begin{defn}
\noindent ~

\noindent $\bullet$ A \emph{projectivity }is a finite product of
perspectivities\emph{.} These mappings are often called \emph{Poncelet
projectivities.}

\noindent \emph{$\bullet$ }When, for example, a projectivity $\pi$
maps the points $A,B,C$ into the points $D,E,F$, in the order written,
we write $ABC\langle\pi\rangle DEF$. 

\noindent \emph{$\bullet$ }We write $\pi\neq\iota$, where $\iota$
is the identity, when, for example, there exists a point $A$ in the
range such that $A^{\pi}\neq A$. 
\end{defn}

\begin{thm}
\label{Thm, projty-preserves-harcons}A projectivity preserves harmonic
conjugates. Thus every Poncelet projectivity is a von Staudt projectivity. 
\end{thm}

\begin{proof}
It will suffice to prove that a harmonic set of points in a range
$\overline{r}$ projects onto a harmonic set of lines in a pencil
$P^{*}$, where $P\notin r$. Given points $A,B,C,D$ on $r$, with
$D=h(A,B;C)$, and $a=PA$, $b=PB$, $c=PC$, $d=PD$, it is required
to show that $d=h(a,b;c)$. Since $A\neq B$, it follows from Proposition
\ref{Prop. lines to two points} that $a\neq b$. By cotransitivity,
it suffices to consider the case in which $C\neq A$. 

Suppose that $h(a,b;c)\neq d$, and suppose further that $C\neq B$.
Thus $c\neq b$. Setting $l=c$, we have $l\neq AB$. Select a point
$R\in b$ so that $R\neq B$ and $R\neq P$. Since $R\neq B=AB\cdot b$,
it follows from Axiom C7 that $R\notin AB$. Since $R\neq P=b\cdot c$,
we have $R\notin l$. Thus the elements $l,R$ may be used in Definition
\ref{Defn. HarConj} to construct the harmonic conjugate $h(A,B;C)$;
it follows that $D=AB\cdot RS$, where $S=AP\cdot BQ$, $P=BR\cdot l$,
and $Q=AR\cdot l.$ 

To construct the harmonic conjugate $h(a,b;c)$, we first set $L=Q$;
then $L\in c$. By Lemma \ref{Lm. HarCon.details}, $L\neq P$, and
by Lemma \ref{Lm. C not A. C not B.}, $L\notin r$. Thus we may use
the elements $L,r$ to construct the line $h(a,b;c)$ using Definition
\ref{Defn. HarConj}(dual). It follows that $h(a,b;c)=(a\cdot b)(r\cdot s)$,
where $s=(a\cdot p)(b\cdot q)$, $p=(b\cdot r)L$, and $q=(a\cdot r)L$.
Now $s=(AP\cdot BQ)(BP\cdot AQ)=SR$, and thus $h(a,b;c)=P(AB\cdot RS)=PD=d$,
contradicting the first assumption above, and negating the second;
hence $C=B$. 

It follows that $c=b$. By Lemma \ref{Lm. HC. C =00003D A} we have
$D=B$, so $d=b$. By the dual of the same lemma, $h(a,b;c)=b$; thus
$h(a,b;c)=d$, contradicting the first assumption above. Hence $h(a,b;c)=d$. 
\end{proof}
The existence of projectivities between ranges will be established
in Theorem \ref{Thm. 3 pts p-projectivity} for the general situation
where it is not known, constructively, whether or not the two ranges
coincide, or, if distinct, whether the common point coincides with
one of the points specified to be mapped. We first consider two lemmas
concerning special situations in which some of this information is
available. 
\begin{lem}
\label{Lm. 2 pts > projection}Let $l$ and $m$ be distinct lines,
with common point $A$. If $Q,R$ are distinct points on $l$, and
$Q',R'$ are distinct points on $m$, with all four points distinct
from $A$, then there exists a projection $\rho:\overline{l}\rightarrow\overline{m}$
such that $AQR\langle\rho\rangle AQ'R'$. 
\end{lem}

\begin{proof}
Since $Q'\neq A=l\cdot m$, it follows from Axiom C7 that $Q'\notin l$,
so $Q'\neq Q$, and $QQ'\neq l$; similarly, $RR'\neq l$. Since $R\neq Q=QQ'\cdot l$,
it follows that $R\notin QQ'$, so $RR'\neq QQ'$. Set $S=QQ'\cdot RR'$;
then $R\neq S$. Now $S\neq R=RR'\cdot l$, so $S\notin l$, and by
symmetry, $S\notin m$. Thus we may define $\rho=\rho(S;l,m)$; it
is clear that $AQR\langle\rho\rangle AQ'R'$. 
\end{proof}
\begin{lem}
\label{Lm. 3 pts > p-projectivity}Let $l$ and $m$ be distinct lines,
with common point $O$. If $P,Q,R$ are distinct points on $l$, and
$P',Q',R'$ are distinct points on $m$, with all six points distinct
from $O$, then there exists a projectivity $\pi:\overline{l}\rightarrow\overline{m}$
such that $PQR\langle\pi\rangle P'Q'R'$. 
\end{lem}

\begin{proof}
Since $Q\neq O=l\cdot m$, it follows that $Q\notin m$; thus $Q\neq P'$,
and similarly for all six points. Set $n=P'Q$; thus $l\cdot n=Q$.
Since $Q\neq R=RR'\cdot l$, it follows from Axiom C7 that $Q\notin RR'$,
so $RR'\neq n$. Thus we may define $R_{0}=RR'\cdot n$. Then $Q\neq R_{0}$,
and by symmetry, $P'\neq R_{0}$. By Lemma \ref{Lm. 2 pts > projection},
there exists a projection $\rho_{1}:\overline{l}\rightarrow\overline{n}$
such that $PQR\langle\rho_{1}\rangle P'QR_{0}$. Also, $n\cdot m=P'$,
so by the same lemma there exists a projection $\rho_{2}:\overline{n}\rightarrow\overline{m}$
such that $P'QR_{0}\langle\rho_{2}\rangle P'Q'R'$. Setting $\pi=\rho_{2}\rho_{1}$,
we obtain $PQR\langle\pi\rangle P'Q'R'$. 
\end{proof}
For a constructive proof of Theorem \ref{Thm. 3 pts p-projectivity},
and also for Lemma \ref{Lm. Conic any 2 base points}, which will
be needed for Pascal's Theorem, we require more points on a line than
has been assumed in Axiom C4. Thus we adopt an additional axiom here.
The determination of the exact number of required points remains an
open problem. \\

\noindent \textbf{Axiom E.\label{Axiom-E. Extension}} Extension.
\emph{There exist at least six distinct points lying on any given
line. }\\

The dual statement is easily verified. 
\begin{thm}
\label{Thm. 3 pts p-projectivity}Given any three distinct points
$P,Q,R$ in a range $\overline{l}$, and any three distinct points
$P',Q',R'$ in a range $\overline{m}$, there exists a projectivity
$\pi:\overline{l}\rightarrow\overline{m}$ such that $PQR\langle\pi\rangle P'Q'R'$.
Similar projectivities exist for other pairs of ranges or pencils. 
\end{thm}

\begin{proof}
Select a point $O_{1}\in l$, distinct from the three given points
in the range $\overline{l}$, and a line $l'$ through $O_{1}$, distinct
from $l$. Select a point $O_{2}\in m$, distinct from the three given
points in $\overline{m}$, and a line $m'$ through $O_{2}$, distinct
from both $m$ and $l'$. Set $O_{3}=l'\cdot m'$, select distinct
points $P_{1},Q_{1},R_{1}\in l'$ distinct from both $O_{1}$ and
$O_{3}$, and select distinct points $P_{2},Q_{2},R_{2}\in m'$ distinct
from both $O_{2}$ and $O_{3}$. Now Lemma \ref{Lm. 3 pts > p-projectivity}
constructs projectivities such that $PQR\rightarrow P_{1}Q_{1}R_{1}\rightarrow P_{2}Q_{2}R_{2}\rightarrow P'Q'R'.$ 
\end{proof}
Classically, at most three perspectivities are needed for Theorem
\ref{Thm. 3 pts p-projectivity}. The constructive proof here uses
three applications of Lemma \ref{Lm. 3 pts > p-projectivity}, each
of which requires two projections, for a total of six. The determination
of the minimum number of perspectivities for a constructive proof
is an open problem. 

\begin{defn}
\noindent A projectivity $\pi$, between two ranges or two pencils,
is said to be \emph{nonperspective} if $x^{\pi}\neq x$ for every
element $x$ in the domain. 
\end{defn}

The following lemma will be needed for Corollary \ref{Cor. projty nonperspective},
and at several places in the study of conics, which will be constructed
using nonperspective projectivities. 
\begin{lem}
\label{Lm. R not S. nonperspective}Let $l$ and $m$ be distinct
lines with common point $O$, let $A,B,C$ be distinct points on $l$,
let $A',B',C'$ be distinct points on $m$, with all six points distinct
from $O$, and define

\[
R=AA'\cdot BB',\,\,\,\,S=BB'\cdot CC',\,\,\,\,n=A'C,
\]
\[
\rho_{1}=\rho(R;l,n),\,\,\,\,\rho_{2}=\rho(S;n,m),\,\,\,\,\pi=\rho_{2}\rho_{1}.
\]

\noindent Then $ABC\langle\pi\rangle A'B'C'$, and the following conditions
are equivalent:

\emph{(a) }The projectivity $\pi:\overline{l}\rightarrow\overline{m}$
is nonperspective. 

\emph{(b)} $O^{\pi}\neq O$.

\emph{(c)} The lines \textup{$AA',BB',CC'$ }\textup{\emph{are nonconcurrent. }}

\emph{(d)} $R\neq S$. 
\end{lem}

\begin{proof}
Except for a change in notation, $\pi$ is the projectivity constructed
in Lemma \ref{Lm. 3 pts > p-projectivity}; thus the definitions are
valid, and $ABC\langle\pi\rangle A'B'C'$. That (a) implies (b) requires
no proof. The equivalence of (c) and (d) is the dual of Proposition
\ref{Prop. lines to two points}. Setting $O_{1}=O^{\rho_{1}}$, we
have $O^{\pi}=O_{1}^{\rho_{2}}=SO_{1}\cdot m$. 

\emph{(b) implies (a).} Given that $O^{\pi}\neq O$, consider any
point $X$ in the range $\overline{l}$. By cotransitivity, either
$X^{\pi}\neq O$ or $X^{\pi}\neq O{}^{\pi}$. In the first case, $X^{\pi}\neq l\cdot m$;
by Axiom C7 it follows that $X^{\pi}\notin l$, and $X^{\pi}\neq X$.
In the second case, $X\neq O=l\cdot m$, so $X\notin m$, and $X\neq X^{\pi}$. 

\emph{(b) implies (d).} Since $O^{\pi}\neq O=RO\cdot m$, it follows
that $O^{\pi}\notin RO$, so $SO_{1}\neq RO=RO_{1}$. Since $R\neq O_{1}=RO_{1}\cdot SO_{1}$,
we have $R\notin SO_{1}$, and $R\neq S$. 

\emph{(d) implies (b).} Given that $R\neq S$, with $R,S\in BB'$,
we have $BB'=RS$; thus $B\in RS$. Since $B\neq O=RO\cdot l$, it
follows that $B\notin RO=RO_{1}$, so $RS\neq RO_{1}$. Since $S\neq R=RS\cdot RO_{1}$,
we have $S\notin RO_{1}$; thus $SO_{1}\neq RO_{1}$. Since $A\neq O$,
we have $A'=A^{\rho_{1}}\neq O^{\rho_{1}}=O_{1}$. Since $O_{1}\neq A'=m\cdot n$,
it follows that $O_{1}\notin m$, and thus $O_{1}\neq O^{\pi}$. Since
$O^{\pi}\neq O_{1}=RO_{1}\cdot SO_{1}$, we have $O^{\pi}\notin RO_{1}$,
and hence $O^{\pi}\neq O$. 
\end{proof}

\section{{\large{}The Fundamental Theorem} \label{SECTION 6. Fund Thm}}

\noindent The \emph{Fundamental Theorem of Projective Geometry} is
the basis for many results, including Pascal's Theorem, the goal of
the present work. The crucial component must be derived from an axiom.
\\

\noindent \textbf{Axiom T.\label{Axiom T. }} \emph{If a projectivity
of a range or pencil onto itself has three distinct fixed elements,
it is the identity. }\\

Classically, Axiom T has the equivalent form, \emph{If a projectivity
$\pi$ from a range onto itself has distinct fixed points $M$ and
$N$, with $\pi\neq\iota$, and $Q$ is a point of the range distinct
from both $M$ and $N$, then $Q^{\pi}\neq Q$.} Constructively, this
appears to be a stronger statement, since the implication ``$\neg(Q^{\pi}=Q)$
implies $Q^{\pi}\neq Q$'' is constructively invalid. To give a proof
of the stronger statement, or to give a Brouwerian counterexample
using the analytic model $\mathbb{P}^{2}(\mathbb{R})$, remains an
open problem. 

\begin{thm}
\emph{\label{Fundamental Theorem}Fundamental Theorem.} Given any
three distinct points $P$, $Q$, $R$ in a range $\overline{l}$,
and any three distinct points $P',Q',R'$ in a range $\overline{m}$,
there exists a unique projectivity $\pi:\overline{l}\rightarrow\overline{m}$
such that $PQR\langle\pi\rangle P'Q'R'$. Similar properties hold
also for other types of projectivity. 
\end{thm}

\begin{proof}
The required projectivity was constructed in Theorem \ref{Thm. 3 pts p-projectivity};
uniqueness follows from Axiom T. 
\end{proof}
\begin{cor}
\label{Cor. Projty-FxPt-perspectivity}If a projectivity from a range
to a distinct range has a fixed point, then it is a perspectivity\emph{.} 
\end{cor}

\begin{proof}
\noindent If $A$ is a fixed point of the projectivity $\pi$, then
it is the point common to the two ranges. Choose distinct points $Q,R$
in the domain, distinct from $A$, denoting the images $Q',R'$; thus
$AQR\langle\pi\rangle AQ'R'$. Use Lemma \ref{Lm. 2 pts > projection}
to construct a projection $\rho$ that agrees with $\pi$ at these
three distinct points; it follows that $\pi=\rho$. 
\end{proof}
\begin{cor}
\label{Cor. projty nonperspective}Let $l$ and $m$ be distinct lines
with common point $O$, let $A,B,C$ be distinct points on $l$, let
$A',B',C'$ be distinct points on $m$, with all six points distinct
from $O$, and let $\pi$ be the projectivity from \textup{$\overline{l}$}
to $\overline{m}$ such that $ABC\langle\pi\rangle A'B'C'$. Then
$\pi$ is nonperspective if and only if the lines \textup{$AA'$,}
$BB'$, \textup{$CC'$} are nonconcurrent. 
\end{cor}

\begin{proof}
Lemma \ref{Lm. R not S. nonperspective} constructs a projectivity
that agrees with $\pi$ at three distinct points. 
\end{proof}
\begin{defn}
\label{Defn. Axis Homology}Let $\pi:\overline{l}\rightarrow\overline{m}$
be a nonperspective projectivity between distinct ranges $\overline{l}$
and $\overline{m}$. Set $O=l\cdot m$, $V=O^{\pi}$, and $U=O^{\pi^{-1}}$.
Pending verification below, the line $UV$ will be called the \emph{axis
of homology for} $\pi$. 
\end{defn}

\begin{thm}
\label{Thm. Axis of homology}Let $l$ and $m$ be distinct lines
with common point $O$, and let $\pi:\overline{l}\rightarrow\overline{m}$
be a nonperspective projectivity. If $A$ and $B$ are distinct points
on $l$, each distinct from $O$, then $A\neq B^{\pi}$, $B\neq A^{\pi}$,
$AB^{\pi}\neq BA^{\pi}$, and the point $AB^{\pi}\cdot BA^{\pi}$
lies on the axis of homology for $\pi$. 
\end{thm}

\begin{proof}
We use the notation of Definition \ref{Defn. Axis Homology}. Since
$\pi$ is nonperspective, we have $V\neq O=l\cdot m$, so it follows
from Axiom C7 that $V\notin l$ and $V\neq U$. Thus the definition
of the axis of homology $UV$ is valid. Since $A\neq O=l\cdot m$,
we have $A\notin m$, so $A\neq B^{\pi}$; similarly, $B\neq A^{\pi}$.
By cotransitivity, either $A\neq U$ or $B\neq U$. In the first case,
we have $A^{\pi}\neq O=l\cdot m$, so $A^{\pi}\notin l$, $OA^{\pi}=m$,
and $BA^{\pi}\neq l$. Since $A\neq B=BA^{\pi}\cdot l$, it follows
that $A\notin BA^{\pi}$, so $AB^{\pi}\neq BA^{\pi}$. By symmetry,
we obtain the same result in the second case. Thus we may set $Q=AB^{\pi}\cdot BA^{\pi}$.

Since $A^{\pi}\notin l$, we have a projection $\rho_{1}:\overline{l}\rightarrow A^{\pi*}$.
Since $A\neq O$, it follows that $A\notin m$, so we also have the
section $\rho_{2}:A^{*}\rightarrow\overline{m}$. Define $\pi_{1}=\rho_{1}\pi^{-1}\rho_{2}$,
a projectivity from $A^{*}$ to $A^{\pi*}$. Since $AA^{\pi}$ is
a fixed line, Corollary \ref{Cor. Projty-FxPt-perspectivity}(dual)
shows that $\pi_{1}$ is a perspectivity. 

Now $(AU)^{\pi_{1}}=A^{\pi}U$, $(AV)^{\pi_{1}}=m$, and $(AB^{\pi})^{\pi_{1}}=A^{\pi}B$.
Thus the axis of the perspectivity $\pi_{1}$ passes through all three
of the points $AU\cdot A^{\pi}U=U$, $AV\cdot m=V$, and $Q$. Hence
$Q$ lies on $UV$. 
\end{proof}
\begin{cor}
\label{Cor. Calculate B pi with axis of homology}Let $l$ and $m$
be distinct lines with common point $O$, let $\pi:\overline{l}\rightarrow\overline{m}$
be a nonperspective projectivity, let $h$ be the axis of homology
for $\pi$, and let $A$ be a point on $l$ with $A\neq O$ and $A^{\pi}\neq O$.
If $B$ is a point on $l$ with $B\neq O$ and $B\neq A$, then $B^{\pi}=A(BA^{\pi}\cdot h)\cdot m$. 
\end{cor}

\begin{proof}
We have $h=UV$, where $V=O^{\pi}$ and $U=O^{\pi^{-1}}$, according
to Definition \ref{Defn. Axis Homology}. Since $\pi$ is nonperspective,
$O^{\pi}\neq O$. Applying $\pi^{-1}$, we have $O\neq O^{\pi^{-1}}=U.$
Since $U\neq O=l\cdot m$, it follows from Axiom C7 that $U\notin m$,
and thus $m\neq h$. Similarly, from $V\neq O$ we obtain $l\neq h$.
From $A^{\pi}\neq O$, applying $\pi^{-1}$, we have $A\neq U=l\cdot h$;
thus $A\notin h$. From $A\neq O$, it follows that $A^{\pi}\neq O^{\pi}=V=m\cdot h$;
thus $A^{\pi}\notin h$, so $BA^{\pi}\neq h$. Now we may set $E=BA^{\pi}\cdot h$;
since $A\notin h$, we have $A\neq E$. By the theorem, $E\in AB^{\pi}$,
so $B^{\pi}\in AE$. Thus $B^{\pi}=AE\cdot m=A(BA^{\pi}\cdot h)\cdot m$. 
\end{proof}
\noindent \emph{Notes for Theorem \ref{Thm. Axis of homology} and
Corollary \ref{Cor. Calculate B pi with axis of homology}.} The restrictions
on the points, or some such, are necessary. For example, if $A=O$
and $B=U$, then the expression $AB^{\pi}$ in the theorem takes the
meaningless form $OO$. In the corollary, if the conditions are all
satisfied except that $A^{\pi}=O$, then the expression $A(BA^{\pi}\cdot h)$
reduces to $UU$. \\

The concept of projectivity may be extended to the entire plane. A
\emph{collineation }of the plane $\mathbb{P}$ is a bijection of the
family $\mathscr{P}$ of points, onto itself, that preserves collinearity
and noncollinearity. A collineation $\sigma$ induces an analogous
bijection $\sigma'$ of the family $\mathscr{L}$ of lines. A collineation
is \emph{projective} if it induces a projectivity on every range and
pencil. 
\begin{prop}
A projective collineation with four distinct fixed points, each three
of which are noncollinear, is the identity. 
\end{prop}

\begin{proof}
Let the collineation $\sigma$ have the fixed points $P,Q,R,S$ as
specified; thus the three distinct lines $PQ,PR,PS$ are fixed. The
mapping $\sigma'$ induces a projectivity on the pencil $P^{*}$;
by the \emph{Fundamental Theorem} this projectivity is the identity.
Thus every line through $P$ is fixed under $\sigma'$; similarly,
the same is true for the other three points. 

Now let $X$ be any point on the plane $\mathbb{P}$. By three successive
applications of cotransitivity for points, we may assume that $X$
is distinct from each of the points $P,Q,R$. Since $PQ\neq PR$ by
Proposition \ref{Prop. lines to two points}, using cotransitivity
for lines we may assume that $XP\neq PQ$. Since $Q\neq P=XP\cdot PQ$,
it follows from Axiom C7 that $Q\notin XP$; thus $XP\neq XQ$. Since
$X=XP\cdot XQ$, and the lines $XP$ and $XQ$ are fixed under $\sigma'$,
it follows that $\sigma X=X$. 
\end{proof}
The construction of a collineation of the plane mapping any set of
four distinct points, each three of which are noncollinear, onto any
similar set of corresponding points, using a constructive synthetic
theory, is an open problem.

\section{{\large{}Involutions} \label{SECTION 7. Involutions}}

\noindent The property of the harmonic conjugate construction, that
the process applied again to the resulting point produces the original
point, Lemma \ref{Lm. hc of D is C}(b), admits a generalization. 
\begin{defn}
An \emph{involution }is a projectivity, from a range or pencil to
itself, of order 2. 
\end{defn}

\begin{thm}
\label{Thm.harcon. involution}Let $A$ and $B$ be distinct points
in a range $\overline{l}$, and let $\upsilon$ be the mapping of
harmonic conjugacy with respect to the base points $A,B$; i.e., set
$X^{\upsilon}=h(A,B;X)$, for all points $X$ in the range $\overline{l}$.
Then $\upsilon$ is an involution. 
\end{thm}

\begin{proof}
Theorem \ref{Thm. harcon.bij} shows that $\upsilon$ is a bijection
of the range $\overline{l}$ onto itself, of order 2. 

To show that $\upsilon$ is a projectivity, we use Definition \ref{Defn. HarConj}
and the notation adopted there. Select a point $R$ outside $AB$,
and select a point $P$ on $BR$ distinct from both $B$ and $R$.
Construct the perspectivities $\rho_{1}(P;AB,AR)$, $\rho_{2}(B;AR,AP)$,
$\rho_{3}(R;AP,AB)$, and the projectivity $\pi_{B}=\rho_{3}\rho_{2}\rho_{1}$. 

Let $X$ be any point in the range $\overline{l}$, with $X\neq B$.
Since $P\neq B=BR\cdot AB$, it follows from Axiom C7 that $P\notin AB$,
so $PX\neq AB$; set $l_{X}=PX$. Since $B\neq X=AB\cdot PX$, we
have $B\notin PX$; thus $BR\neq PX$. Now $R\neq P=BR\cdot PX$,
so $R\notin l_{X}$. Thus, to construct the harmonic conjugate $X^{\upsilon}$,
we may use the point $R$ and the line $l_{X}$ in the definition.
We find that $X^{\pi_{B}}$ = $R(B(PX\cdot AR)\cdot AP)\cdot AB$
= $R(B(l_{X}\cdot AR)\cdot AP)\cdot AB$ = $R(BQ_{X}\cdot AP)\cdot AB$
= $RS_{X}\cdot AB$ = $X^{\upsilon}$. Thus $\pi_{B}$ agrees with
$\upsilon$ for all points in the range $\overline{l}$ that are distinct
from $B$. 

Similarly, construct the projectivity $\pi_{A}$; it will agree with
$\upsilon$ for all points in $\overline{l}$ that are distinct from
$A$. Choose any three distinct points on $l$, each distinct from
both $A$ and $B$. Since $\pi_{A}$ and $\pi_{B}$ agree at these
three points, by the \emph{Fundamental Theorem} they are the same
projectivity; call it $\pi$. By cotransitivity, each point in the
range $\overline{l}$ is either distinct from $A$, or distinct from
$B$. Hence $\upsilon=\pi$. 
\end{proof}
\begin{lem}
\label{Lm. 4 pts - projectivity}Given any four distinct points $A,B,C,D$
in a range $\overline{l}$, there exists a projectivity $\pi$ from
$\overline{l}$ to itself such that $ABCD\langle\pi\rangle BADC$. 
\end{lem}

\begin{proof}
Select a line $m$ through $D$, distinct from $l$, and select a
point $Q$ outside both lines $l$ and $m$. Set $\rho_{1}=\rho(Q;l,m)$.
Since $A\neq D=l\cdot m$, it follows from Axiom C7 that $A\notin m$,
and thus $AQ\neq m$. Similarly, both $B$ and $C$ lie outside $m$,
and $BQ$ and $CQ$ are both distinct from $m$. Since $A\neq C=CQ\cdot l$,
we have $A\notin CQ$. Similarly, $B\notin CQ$, so $BQ\neq CQ$.
Set $\rho_{2}=\rho(A;m,CQ)$, $R=AQ\cdot m$, $S=BQ\cdot m$, and
$T=CQ\cdot m$. Since $S\neq Q=BQ\cdot CQ$, it follows that $S\notin CQ$
and $AS\neq CQ$; set $U=AS\cdot CQ$. Since $D\neq B=BQ\cdot l$,
we have $D\notin BQ$, so $D\neq S$. Since $S\neq D=l\cdot m$, we
have $S\notin l$, so $S\neq A$. Finally, set $\rho_{3}=\rho(S;CQ,l)$.
It is clear that $ABCD\langle\rho_{1}\rangle RSTD\langle\rho_{2}\rangle QUTC\langle\rho_{3}\rangle BADC$;
set $\pi=\rho_{3}\rho_{2}\rho_{1}$. 
\end{proof}

\begin{thm}
\label{Thm. Interchange 2 - involution}A projectivity from a range
to itself, that interchanges two distinct elements, is an involution. 
\end{thm}

\begin{proof}
Let $AB\langle\pi\rangle BA$, where $A\neq B$. Let $X$ be any point
in the range, and set $Y=X^{\pi}$. By symmetry and cotransitivity,
it suffices to consider the case in which $X\neq A$. Suppose that
$Y^{\pi}\ne X$, and suppose further that $Y\neq X$ and $X\neq B$. 

Applying $\pi$, it follows that $Y\neq A$ and $Y\neq B$. Using
Lemma \ref{Lm. 4 pts - projectivity}, construct a projectivity $\pi_{1}$
with $ABXY\langle\pi_{1}\rangle BAYX$. Since the projectivity $\pi$
agrees with $\pi_{1}$ at three distinct points, by the \emph{Fundamental
Theorem} we have $\pi=\pi_{1}$, so $Y^{\pi}=X$. This contradicts
the first assumption above, and negates the last. 

Thus $X=B$; applying $\pi$ here, we have $Y=A$. Thus $Y^{\pi}=B=X$,
contradicting the first assumption and negating the second. Thus $Y=X$,
and it follows that $Y^{\pi}=X^{\pi}=Y=X$, contradicting and negating
the first assumption. Hence $Y^{\pi}=X$. 
\end{proof}
\begin{thm}
If an involution $\pi$ on a range $\overline{l}$ has a fixed point
$M$, then it has a second, distinct, fixed point $N$, and $\pi$
is the mapping of harmonic conjugacy with respect to these points;
thus $X^{\pi}=h(M,N;X)$, for all points $X$ in $\overline{l}$. 
\end{thm}

\begin{proof}
Select a point $A$ in $\overline{l}$ so that $A^{\pi}\neq A$. Either
$M\neq A$ or $M\neq A^{\pi}$; in either case, applying $\pi$ we
find that $M$ is distinct from both $A$ and $A^{\pi}$. Set $N=h(A,A^{\pi};M)$;
by Lemma \ref{Lm. C not ei base point}, $N\neq M$, and by Lemma
\ref{Lm. hc of D is C}(a), we have $N=h(A^{\pi},A;M)$. Applying
the projectivity $\pi$, we have $N^{\pi}=h(A,A^{\pi};M)=N$. Thus
$N$ is a second fixed point.

If an alternative selection of the point $A$ results in the second
fixed point $N_{1}\neq M$, then, by the \emph{Fundamental Theorem,}
$N_{1}=N$, since $N_{1}\neq N$ would mean that we have three distinct
fixed points. Thus, for any point $X$ in the range $\overline{l}$,
if $X^{\pi}\neq X$, then $h(X,X^{\pi};M,N)$. 

Now consider any point $X$ in the range $\overline{l}$. By cotransitivity,
we may assume that $X\neq M$. Suppose that $X^{\pi}\neq h(M,N;X)$,
and suppose further that $X\neq N$ and $X\neq X^{\pi}$. We have
$h(X,X^{\pi};M,N)$ from above. By Theorem \ref{Thm, projty-preserves-harcons}
and Lemma \ref{Lm. 4 pts - projectivity}, it follows that $h(M,N;X,X^{\pi})$,
contradicting the first assumption, and negating the last; thus $X=X^{\pi}$.
Now we have three distinct fixed points; this contradiction negates
the second assumption. Thus $X=N$, a contradiction negating the first
assumption. Hence $X^{\pi}=h(M,N;X)$. 
\end{proof}

\section{{\large{}Conics} \label{SECTION 8. Conics}}

\noindent We define conics by means of projectivities, using the method
of Steiner {[}Ste32{]}. 
\begin{defn}
\label{Defn. Conic}Let $\pi:U^{*}\rightarrow V^{*}$ be a nonperspective
projectivity between distinct pencils of lines. The \emph{conic $\kappa=\kappa(\pi;U,V)$
defined by $\pi$ }is the locus of points $\{l\cdot l^{\pi}:l\in U^{*}\}$.\footnote{With $\pi$ being nonperspective, this is usually called a \emph{non-singular
conic}; we leave the constructive study of the singular conics for
a later time. } For any point $X$, we say that \emph{$X$ lies outside $\kappa$},
written $X\notin\kappa$, if $X\neq Y$ for all points $Y$ on $\kappa$.
At times, the locus $\kappa$ may be called a \emph{point-conic};
the dual locus is a \emph{line-conic}. 
\end{defn}

\begin{prop}
\label{Prop. conic 4  parts}Let \emph{$\kappa=\kappa(\pi;U,V)$}\textup{
}\textup{\emph{be a conic.}}\textup{ }

(a) The base points $U$ and $V$ are points of $\kappa$.

(b) Any three distinct points on \textup{$\kappa$} are noncollinear. 

(c) For any point $P$ on $\kappa$, the line $l\in U^{*}$, such
that $P=l\cdot l^{\pi}$, is unique. If $P\neq U$, then $l=UP$,
while if $P\neq V$, then $l=(VP)^{\pi^{-1}}$. 

(d) For any point $X$, if $\neg(X\notin\kappa)$, then $X\in\kappa$. 
\end{prop}

\begin{proof}
(a) The line $o=UV$ in the pencil $U^{*}$ has a corresponding line
$o^{\pi}$ in $V^{*}$; this pair of lines determines the point $o\cdot o^{\pi}=V$
of $\kappa$. Similarly for $U$. 

(b) This follows from Corollary \ref{Cor. projty nonperspective}(dual). 

(c) By cotransitivity, either $P\neq U$ or $P\neq V$. In the first
case, both $U$ and $P$ lie on $l$, so $l=UP$. Similarly, in the
second case we have $l^{\pi}=VP$. 

(d) Let $X$ be a point on the plane such that $\neg(X\notin\kappa)$.
By cotransitivity, we may assume that $X\neq U$; set $x=UX$ and
$Z=x\cdot x^{\pi}$. Suppose that $X\neq Z$. We will show that $X\neq Y$
for any point $Y$ of $\kappa$. Either $Y\neq X$ or $Y\neq U$.
We need to consider only the second case; set $y=UY$. By (c), we
have $Y=y\cdot y^{\pi}$. Either $Y\neq X$ or $Y\neq Z$; again we
need to consider only the second case. Since $Y\neq x\cdot x^{\pi}$,
it follows from Axiom C7 that either $Y\notin x$ or $Y\notin x^{\pi}$.
In the first subcase, $y\neq x$. In the second subcase, $y^{\pi}\neq x^{\pi}$,
and since $\pi$ is a bijection we again have $y\neq x$. Since $X\neq U=x\cdot y$,
it follows that $X\notin y$, and $X\neq Y$. This shows that $X\notin\kappa$,
a contradiction. Hence $X=Z$, and now we have $X\in\kappa.$ 
\end{proof}
\begin{prop}
\label{Prop. 5 points determine conic with fixed base points}Given
any five distinct points U,V,A,B,C, each three of which are noncollinear,
there exists a unique conic, with base points U and V, containing
all five points. 
\end{prop}

\begin{proof}
This follows from Corollary \ref{Cor. projty nonperspective}(dual). 
\end{proof}
The following three lemmas are required for Pascal's Theorem. The
first is a special case; two vertices of the hexagon are the base
points of a projectivity that defines the conic. 
\begin{lem}
\label{Conic Lemma 1}Let $\kappa=\kappa(\pi;U,V)$ be any conic,
and $A,B,C,X$ points of $\kappa$, with all six points distinct.
Then the three points $O=UA\cdot VC$, $Y=UX\cdot BC$, $Z=VX\cdot AB$
are distinct and collinear. 
\end{lem}

\begin{proof}
That the three points in question are properly defined follows from
Proposition \ref{Prop. conic 4  parts}(b). Set $D=UA\cdot BC$ and
$E=VC\cdot AB$; thus $D\neq A$, so $DA=UA$, and $E\neq C$, so
$CE=VC$. Consider the sections $\rho_{1}:U^{*}\rightarrow\overline{BC}$
and $\rho_{2}:V^{*}\rightarrow\overline{AB}$; thus $UA,UB,UC,UX\langle\rho_{1}\rangle DBCY$
and $VA,VB,VC,VX\langle\rho_{2}\rangle ABEZ$. Setting $\pi_{1}=\rho_{2}\pi\rho_{1}^{-1}$,
we obtain a projectivity $\pi_{1}:\overline{BC}\rightarrow\overline{AB}$,
with $DBCY\langle\pi_{1}\rangle ABEZ$. Since $B$ is a fixed point,
it follows from Corollary \ref{Cor. Projty-FxPt-perspectivity} that
$\pi_{1}$ is a projection; the center is $DA\cdot CE=UA\cdot VC=O$.
The center of a projection lies outside each of the ranges it maps;
thus the three points are distinct. Also, $Z=Y^{\pi_{1}}=OY\cdot AB$,
and hence $Z\in OY$. 
\end{proof}
\begin{lem}
\label{Conic Lemma 2}Let $U,A,B,C,V,X$ be distinct points, each
three of which are noncollinear. If the three points $O=UA\cdot VC$,
$Y=UX\cdot BC$, $Z=VX\cdot AB$ are collinear, then there exists
a nonperspective projectivity $\pi:U^{*}\rightarrow V^{*}$ such that
all six points are on the conic $\kappa(\pi;U,V)$. 
\end{lem}

\begin{proof}
Set $a=UA$, $b=UB$, $c=UC$, $x=UX$, $a'=VA$, $b'=VB$, $c'=VC$,
$x'=VX$, $l=BC$, $m=AB$, $D=a\cdot l$, and $E=c'\cdot m$; thus
$O=a\cdot c'$, $Y=x\cdot l$, and $Z=x'\cdot m$. 

Since $C\notin a$, we have $C\neq O$. Since $O\neq C=c'\cdot l$,
it follows from Axiom C7 that $O\notin l$. Similarly, $O\notin m$,
so we may construct the projection $\rho=\rho(O;l,m)$. We have $D^{\rho}=OD\cdot m=a\cdot m=A$,
and $C^{\rho}=OC\cdot m=c'\cdot m=E$. By hypothesis, $Z\in OY$,
so $Y^{\rho}=OY\cdot m=Z$. Thus $DBCY\langle\rho\rangle ABEZ$. We
also have sections of the pencils, $\rho_{1}:U^{*}\rightarrow\overline{l}$
and $\rho_{2}:V^{*}\rightarrow\overline{m}$, such that $abcx\langle\rho_{1}\rangle DBCY$
and $a'b'c'x'\langle\rho_{2}\rangle ABEZ$. Setting $\pi=\rho_{2}^{-1}\rho\rho_{1}$,
we obtain a projectivity $\pi:U^{*}\rightarrow V^{*}$, with $abcx\langle\pi\rangle a'b'c'x'$.
Since the points $a\cdot a'$, $b\cdot b'$, $c\cdot c'$ are noncollinear,
it follows from Corollary \ref{Cor. projty nonperspective}(dual)
that $\pi$ is nonperspective. The conic $\kappa(\pi;U,V)$ clearly
includes all six points.
\end{proof}
\begin{lem}
\emph{(Steiner)} \label{Lm. Conic any 2 base points}Given any conic
$\kappa=\kappa(\pi;U,V)$, and any distinct points $U_{1},V_{1}$
on $\kappa$, there exists a unique nonperspective projectivity \textup{$\pi_{1}:U_{1}^{*}\rightarrow V_{1}^{*}$}
such that $\kappa=\kappa(\pi_{1};U_{1},V_{1})$. 
\end{lem}

\begin{proof}
(a) \emph{Special case; $U_{1}=U$ and $V_{1}=W$, where $W$ is any
point of $\kappa$ that is distinct from both $U$ and $V$.} Select
points $A,B$ on $\kappa$ such that the points $U,A,W,B,V$ are distinct,
and let $X$ be any point of $\kappa$ distinct from these five. Applying
Lemma \ref{Conic Lemma 1} to the points $U,A,W,B,V,X$, we obtain
three collinear points of interest. A cyclic permutation of these
six points results in $W,B,V,X,U,A$, with the same three points of
interest. Thus Lemma \ref{Conic Lemma 2} applies, and we obtain a
conic $\kappa_{1}=\kappa(\pi_{1};U,W)$ containing all six points.
By Proposition \ref{Prop. 5 points determine conic with fixed base points},
this conic, with base points $U,W$, is independent of the choice
of $X$. 

Let $Y$ be any point of $\kappa$, and suppose that $Y\notin\kappa_{1}$.
Now $Y$ is distinct from each of the points $U,V,W,A,B$, so the
construction of $\kappa_{1}$ may be repeated with $Y$ in place of
$X$; thus $Y$ is on $\kappa_{1}$, a contradiction. From Proposition
\ref{Prop. conic 4  parts}(d), it follows that $Y\in\kappa_{1}$;
thus $\kappa\subset\kappa_{1}$. 

Applying the same construction method, now to $\kappa_{1}$, we obtain
a conic $\kappa_{2}=\kappa(\pi_{2};U,V)$, with $\kappa_{1}\subset\kappa_{2}$.
By Proposition \ref{Prop. 5 points determine conic with fixed base points},
$\kappa=\kappa_{2}$, and hence $\kappa=\kappa_{1}$. 

(b) \emph{General case. }Using Axiom E, select distinct points $W_{1},W_{2}$
on $\kappa$, each distinct from all four points $U,V,U_{1},V_{1}$.
Four applications of special case (a) result in the sequence $(U,V),$
$(U,W_{2}),$ $(W_{1},W_{2}),$ $(U_{1},W_{2}),$ $(U_{1},V_{1})$
of base-point changes.

(c) The uniqueness follows from the \emph{Fundamental Theorem}. 
\end{proof}
\begin{thm}
There exists a unique conic containing any given five distinct points,
each three of which are noncollinear. 
\end{thm}

\begin{proof}
The existence follows from Theorem \ref{Thm. 3 pts p-projectivity}
and Corollary \ref{Cor. projty nonperspective}(dual), the uniqueness
from Proposition \ref{Prop. 5 points determine conic with fixed base points}
and Lemma \ref{Lm. Conic any 2 base points}.
\end{proof}

\section{{\large{}Pascal's Theorem} \label{SECTION 9. Pascal}}

\noindent For information concerning Blaise Pascal (1623-62), see
{[}Kli72, p. 295-299{]}. 
\begin{defn}
A \emph{simple hexagon} is a set \emph{$ABCDEF$} of six distinct
points, in cyclic order, each three of which are noncollinear. The
six points are the \emph{vertices}; the six lines joining successive
points are the \emph{sides}. The pairs of sides $(AB,DE)$, $(BC,EF)$,
$(CD,FA)$ are said to be \emph{opposite}. 
\end{defn}

\begin{thm}
\emph{(Pascal)} If a simple hexagon is inscribed in a conic, the three
points of intersection of the pairs of opposite sides are distinct
and collinear. 
\end{thm}

\begin{proof}
Label the inscribed hexagon as $UABCVX$, and apply Lemma \ref{Lm. Conic any 2 base points}
to view the conic as $\kappa=\kappa(\pi;U,V)$. Now Lemma \ref{Conic Lemma 1}
yields the result. 
\end{proof}
According to legend, Pascal gave in addition some four hundred corollaries.
Here we have only one; it recalls a traditional construction method
for drawing a conic ``point by point''.\footnote{For example, as in {[}You30, p. 68{]}.}
\begin{cor}
\label{Cor. Pascal 5 points- sixth}Let $A,B,C,D,E$ be five distinct
points of a conic $\kappa$, and let $l$ be a line through $E$ that
avoids $A,B,C,D$. If $l$ passes through a distinct sixth point F
of $\kappa$, then $F=l\cdot A(CD\cdot(AB\cdot DE)(BC\cdot l))$. 
\end{cor}

\begin{proof}
The \emph{Pascal line} $p$ of the hexagon $ABCDEF$ passes through
the three points $X=AB\cdot DE$, $Y=BC\cdot EF$, and $Z=CD\cdot AF$.
Since $A\notin CD$, we have $A\neq Z$, so $AF=AZ$. Since $B\notin CD$,
we have $BC\neq CD$, so by cotransitivity for lines either $p\neq BC$
or $p\neq CD$. In the first case, since $C\notin EF$, we have $C\neq Y=BC\cdot p$,
so it follows from Axiom C7 that $C\notin p$. Thus in both cases
we have $CD\neq p$, and $Z=CD\cdot p$. Now $F=EF\cdot AF=l\cdot AZ=l\cdot A(CD\cdot p)=l\cdot A(CD\cdot XY)=l\cdot A(CD\cdot(AB\cdot DE)(BC\cdot l)$.
\end{proof}

\section{{\large{}Tangents and secants} \label{SECTION 10. Tangents Secants} }

\noindent The construction of poles and polars with respect to a conic,
in Section \ref{SECTION 11. Poles and polars}, will involve the properties
of tangents and secants. 
\begin{defn}
\label{DEFN. tangent}Let $\kappa$ be a conic, and $P$ a point on
$\kappa$. A line $t$ that passes through $P$ is said to be \emph{tangent
to $\kappa$ at $P$} if $P$ is the unique point of $\kappa$ that
lies on $t$. The dual concept is a \emph{point of contact} $L$ of
a line $l$ that belongs to a line-conic $\lambda$. 
\end{defn}

\begin{prop}
\label{Prop.  tangent to conic}

\label{Prop.  tangent to conic-1}Let $\kappa$ be a conic, $P$ a
point on $\kappa$, and $t$ a line passing through $P$. The following
are equivalent.

(a) The line $t$ is tangent to $\kappa$ at $P$. 

(b) For any point $Q$ of $\kappa$, if $Q\neq P$ and $\pi$ is the
nonperspective projectivity such that $\kappa=\kappa(\pi;Q,P)$, then
$t=(QP)^{\pi}$. 

(c) There exists a point $Q$ of $\kappa$ with $Q\neq P$, and corresponding
nonperspective projectivity $\pi$, with $\kappa=\kappa(\pi;Q,P)$,
such that \textup{$t=(QP)^{\pi}$}. 
\end{prop}

\begin{proof}
Given (a), and a point Q as specified in (b), set $u=t^{\pi^{-1}}$.
Since $u\cdot t$ is a point of $\kappa$ which lies on $t$, it must
be $P$; thus $u=QP$, so $t=u^{\pi}=(QP)^{\pi}$. 

Given (c), and a point $Q$ as specified, let $R$ be a point on $\kappa$
with $R\in t$. Suppose that $R\neq P$; thus $t=PR$. From Proposition
\ref{Prop. conic 4  parts}(c), applied to $\pi^{-1}:P^{*}\rightarrow Q^{*}$,
it follows that $R=t\cdot t^{\pi^{-1}}=t\cdot QP=P$, a contradiction;
hence $R=P$, and this proves (a). 
\end{proof}
\begin{cor}
\label{Prop. tangent exists}Let $\kappa$ be a conic, and $P$ any
point on $\kappa$.

(a) There exists a unique line $t$ that is tangent to $\kappa$ at
$P$. 

(b) Let $t$ be the tangent to $\kappa$ at $P$. If $Q$ is any point
on $\kappa$, with $Q\neq P$, then $Q\notin t$. 
\end{cor}

\begin{proof}
(a) follows directly from Lemma \ref{Lm. Conic any 2 base points}
and Proposition \ref{Prop.  tangent to conic}. For (b), select any
point Q as specified in (b) of the same proposition; thus $t=(QP)^{\pi}$.
Since $\pi$ is nonperspective, we have $t\neq QP$. Since $Q\neq P=t\cdot QP$,
it follows from Axiom C7 that $Q\notin t$. 
\end{proof}
\begin{thm}
\label{Thm. 5 pts + 1 tangent}Let $\kappa$ be any conic, and let
$UABCV$ be a pentagon inscribed in $\kappa$, with five distinct
vertices. The point of intersection of the tangent $u$ at $U$, with
the side opposite, is collinear with the points of intersection of
the other two pairs of nonadjacent sides. That is, the three points
$O=UA\cdot VC$, $Z=UV\cdot AB$, $Y=u\cdot BC$ are collinear. 
\end{thm}

\begin{proof}
Using Lemma \ref{Lm. Conic any 2 base points}, construct the projectivity
$\pi$ so that $\kappa=\kappa(\pi;U,V)$. That the three points in
question are properly defined follows from Proposition \ref{Prop. conic 4  parts}(b)
and Proposition \ref{Prop. tangent exists}(b). Set $D=UA\cdot BC$,
and $E=VC\cdot AB$. Since $A\notin BC$, we have $A\neq D$, so $DA=UA$.
Since $C\notin AB$, we have $C\neq E$, so $CE=VC$. 

Construct the sections $\rho_{1}:U^{*}\rightarrow\overline{BC}$ and
$\rho_{2}:V^{*}\rightarrow\overline{AB}$; clearly, $UA,UB,UC,u$
$\langle\rho_{1}\rangle DBCY$ and $VA,VB,VC,VU\langle\rho_{2}\rangle ABEZ$.
Setting $\pi_{1}=\rho_{2}\pi\rho_{1}^{-1}$, we obtain a projectivity
$\pi_{1}:\overline{BC}\rightarrow\overline{AB}$, with $DBCY\langle\pi_{1}\rangle ABEZ$.
Since $B$ is a fixed point, it follows from Corollary \ref{Cor. Projty-FxPt-perspectivity}
that $\pi_{1}$ is a projection; the center is $DA\cdot CE=UA\cdot VC=O$.
Thus $Z=Y^{\pi_{1}}\in OY$. 
\end{proof}
\begin{thm}
\label{Thm. 4  pts + 2  tangents}Let $\kappa$ be any conic, and
let $UABV$ be a quadrangle inscribed in $\kappa$. The point of intersection
of the tangent at $U$ with the side $VB$, the point of intersection
of the tangent at $V$ with the side $UA$, and the diagonal point
lying on $UV$, are collinear. 
\end{thm}

\begin{proof}
Let $\pi$ be the projectivity such that $\kappa=\kappa(\pi;U,V)$,
and denote the tangents at $U$ and $V$ by $u$ and $v$. The three
points in question are then $Y=u\cdot VB$, $O=v\cdot UA$, and $Z=UV\cdot AB$.
Set $D=UA\cdot VB$, and $E=v\cdot AB$. Since $A\notin VB$, we have
$A\neq D$, so $DA=UA$. Since $V\notin AB$, we have $V\neq E$,
so $VE=v$.

Construct the sections $\rho_{1}:U^{*}\rightarrow\overline{VB}$ and
$\rho_{2}:V^{*}\rightarrow\overline{AB}$; it is clear that $UA,UB,UV,u$\emph{
$\langle\rho_{1}\rangle DBVY$} and $VA,VB,v,VU\langle\rho_{2}\rangle ABEZ$.
Setting $\pi_{1}=\rho_{2}\pi\rho_{1}^{-1}$, we obtain a projectivity
$\pi_{1}:\overline{VB}\rightarrow\overline{AB}$, with $DBVY\langle\pi_{1}\rangle ABEZ$.
Since $B$ is a fixed point, $\pi_{1}$ is a projection with center
$DA\cdot VE=UA\cdot v=O$. Thus $Z=Y^{\pi_{1}}\in OY$. 
\end{proof}
\begin{thm}
\label{Thm. 4 pts - tangents intersection point}Let $\kappa$ be
any conic, and let $UABV$ be a quadrangle inscribed in $\kappa$.
The point of intersection of the tangents at $U$ and $V$, and the
two diagonal points not lying on $UV$, are distinct and collinear. 
\end{thm}

\begin{proof}
Construct the projectivity $\pi$ so that $\kappa=\kappa(\pi;U,V)$.
Denote the tangents at $U$ and $V$ by $u$ and $v$; the three points
in question are then $O=u\cdot v$, $D_{1}=UA\cdot VB$, $D_{2}=UB\cdot VA$.
Set $E=u\cdot VB$ and $F=v\cdot UB$. Since $U\notin VB$, we have
$U\neq E$, so $EU=u$; by symmetry, $FV=v$. 

Construct the sections $\rho_{1}:U^{*}\rightarrow\overline{VB}$ and
$\rho_{2}:V^{*}\rightarrow\overline{UB}$; clearly, $u,UB,UV,UA$
$\langle\rho_{1}\rangle EBVD_{1}$ and $VU,VB,v,VA\langle\rho_{2}\rangle UBFD_{2}$.
Setting $\pi_{1}=\rho_{2}\pi\rho_{1}^{-1}$, we obtain a projectivity
$\pi_{1}:\overline{VB}\rightarrow\overline{UB}$, with $EBVD_{1}\langle\pi_{1}\rangle UBFD_{2}$.
Since $B$ is a fixed point, $\pi_{1}$ is a projection; the center
is $EU\cdot FV=u\cdot v=O$. The center of a projection lies outside
each of the ranges it maps; thus the three points are distinct. Since
$D_{2}=D_{1}^{\pi_{1}}\in OD_{1}$, the points are collinear.
\end{proof}
\begin{cor}
Given any conic $\kappa=\kappa(\pi;U,V)$, the center of homology
of the projectivity $\pi$ is the intersection $u\cdot v$ of the
tangents to $\kappa$ at $U$ and $V$. 
\end{cor}

The following theorem is related to the existence of secants, to be
constructed in Theorem \ref{Thm. 2 secants}.
\begin{thm}
\label{THM. tangents and line-conic}The following two statements
are equivalent:

(a) The tangents at any three distinct points of a point-conic are
nonconcurrent; the points of contact of any three distinct lines of
a line-conic are noncollinear. 

(b) The family of all tangents to a point-conic is a line-conic; the
family of all points of contact of a line-conic is a point-conic. 
\end{thm}

\begin{proof}
Since (a) follows directly from (b) and Proposition \ref{Prop. conic 4  parts}(b),
with its dual, it only remains to prove that (a) implies (b).

Given a point-conic $\kappa=\kappa(\pi;A,B)$, select a point $C\in\kappa$,
distinct from both $A$ and $B$, and let $a,b,c$ be the tangents
at these three points. Since $A\neq B$, by Proposition \ref{Prop. tangent exists}(b)
we have $A\notin b$, and thus $a\neq b$; similarly for the other
points and tangents. Set $E=a\cdot b$, $F=b\cdot c$, and $G=a\cdot c$.
It follows from (a) that the points $A,E,G$ are distinct, as are
$E,B,F$. Thus we may construct the projectivity $\varphi:\overline{a}\rightarrow\overline{b}$
such that $AEG\langle\varphi\rangle EBF$. Since $E^{\varphi}\neq E$,
it follows from Lemma \ref{Lm. R not S. nonperspective} that $\varphi$
is nonperspective, so the family of lines $\lambda=\lambda(\varphi;a,b)=\{QQ^{\varphi}:Q\in a\}$
is a line-conic; the axis of homology for $\varphi$ is clearly $h=AB$. 

(1) \emph{If $P$ is any point of the point-conic $\kappa$ with $P\neq A,B,C$,
then the tangent $p$ to $\kappa$ at $P$ is a line of the line-conic
$\lambda$. }To prove this, denote the diagonals of the quadrangle
$ABCP$ by $D_{1}=AC\cdot BP$, $D_{2}=AB\cdot CP$, and $D_{3}=AP\cdot BC$.
Also, set $S=a\cdot p$ and $T=b\cdot p$. By Theorem \ref{Thm. 4 pts - tangents intersection point},
the points $F,D_{1},D_{2}$ are distinct and collinear, as are the
points $S,D_{1},D_{2}$, the points $G,D_{1},D_{2}$ and the points
$T,D_{1},D_{2}$. 

Since $F\neq E=a\cdot b$, it follows from Axiom C7 that $F\notin a$,
so $F\neq S$ and $SF=SD_{2}$. From (a), we have $G\neq E$, $F\neq E$,
$S\neq E$, and $S\neq G$, so we may apply Corollary \ref{Cor. Calculate B pi with axis of homology}
and the axis of homology. Thus $S^{\varphi}=G(SF\cdot h)\cdot b=G(SD_{2}\cdot AB)\cdot b=GD_{2}\cdot b=T$,
and $SS^{\varphi}=ST=(a\cdot p)(b\cdot p)=p$; hence $p\in\lambda$. 

(2) \emph{The tangents $a,b,c$ of the point-conic $\kappa$ are each
lines of the line-conic $\lambda$. }For $a$ and $b$, this follows
from Proposition \ref{Prop. conic 4  parts}(a)(dual). For $c$, it
suffices to note that $GF=c$. 

(3) \emph{Each tangent to the point-conic $\kappa$ is a line of the
line-conic $\lambda$. }Let $P$ be any point of $\kappa$, with tangent
$p$, and suppose that $p\notin\lambda$. Suppose further, in succession,
that $P\neq A$, $P\neq B$, $P\neq C$. Now we have $p\in\lambda$
by (1), a contradiction. Thus $P=C$, so $p=c$, and $p\in\lambda$
by (2), a contradiction. Thus $P=B$. Continuing this way, we arrive
at a final contradiction. Thus $\neg(p\notin\lambda)$, and it follows
from Proposition \ref{Prop. conic 4  parts}(d)(dual) that $p\in\lambda$. 

(4) \emph{Each point of contact of $\lambda$ is a point of $\kappa$.}
Apply the dual of the method in (1) to the line-conic $\lambda$,
using the second part of condition (a). The result is a point-conic
$\kappa_{1}=\kappa(\psi:A,B)$, where $aeg\langle\psi\rangle ebf$,
and $e=AB$, $f=BC$, $g=AC$. It follows from the dual of (3) that
every point of contact of $\lambda$ is a point of $\kappa_{1}$.
We have $a,AB,AC\langle\pi\rangle BA,b,BC$, and these six lines,
in order, are identical to those just noted for $\psi$. Thus, by
the \emph{Fundamental Theorem}, $\psi=\pi$ and $\kappa_{1}=\kappa$. 

(5) \emph{Each line of the line-conic $\lambda$ is a tangent of the
point-conic $\kappa$.} Let $l$ be a line of $\lambda$, with point
of contact $L$. By Definition \ref{DEFN. tangent}, $l$ is the unique
line of $\lambda$ passing through $L$, and by (4), $L\in\kappa$.
Let $t$ denote the tangent to $\kappa$ at $L$. By (3), $t$ is
a line of $\lambda$; hence $t=l$. 
\end{proof}
A line that passes through two distinct points of a conic $\kappa$
is a \emph{secant} \emph{of} $\kappa$. 
\begin{lem}
\label{Lm. second point}Let $\kappa$ be a conic, $P$ a point on
$\kappa$, and $t$ the tangent to $\kappa$ at $P$. If $l$ is a
line through $P$, and $l\neq t$, then $l$ passes through a second
point $R$ of $\kappa$, distinct from $P$; thus $l$ is a secant
of $\kappa$. 
\end{lem}

\begin{proof}
Select $Q$ and $\pi$ as in Proposition \ref{Prop.  tangent to conic}(c).
Thus $\kappa=\kappa(\pi;Q,P)$ and $t=(QP)^{\pi}$, so $P\in t^{\pi^{-1}}$.
Set $R=l\cdot l^{\pi^{-1}}$. Since $l\neq t$ and $P\neq Q=l^{\pi^{-1}}\cdot t^{\pi^{-1}}$,
it follows from Axiom C7 that $P\notin l^{\pi^{-1}}$, and hence $P\neq R$. 
\end{proof}
For any conic, the following theorem will provide, through an arbitrary
point of the plane, the one secant needed to construct polars in Lemma
\ref{Thm. Construct polar}, and the two distinct secants needed for
Corollary \ref{Cor. polar from quadrangle}, relating polars to inscribed
quadrangles. The need for this theorem contrasts with complex geometry,
where every line meets every conic.
\begin{thm}
\label{Thm. 2 secants}Let $\kappa$ be a conic, and assume statement
(a) of Theorem \ref{THM. tangents and line-conic}. 

(a) Through any given point of the plane, we may construct at least
two distinct secants of $\kappa$. 

(b) On any given line of the plane, we may construct at least two
distinct points, through each of which pass two tangents of $\kappa.$
\end{thm}

\begin{proof}
(a) Let $P$ and $\kappa$ be given, and select distinct points $A,B,C$
on $\kappa$, with tangents $a,b,c$. By hypothesis, these tangents
are nonconcurrent; thus the points $E=a\cdot b$ and $F=b\cdot c$
are distinct. Either $P\neq E$ or $P\neq F$; it suffices to consider
the first case. Then, by Axiom C7, either $P\notin a$ or $P\notin b$.
It suffices to consider the first subcase; thus $P\neq A$ and $PA\neq a$.
Now it follows from Lemma \ref{Lm. second point} that $PA$ is a
secant. 

Denote the second point of $PA$ that lies on $\kappa$ by $R$, and
choose distinct points $A',B',C'$ on $\kappa$, each distinct from
both $A$ and $R$. Using these three points, construct a secant through
$P$ with the above method; we may assume that it is $PA'$. By Proposition
\ref{Prop. conic 4  parts}(b), $A'\notin AR=PA$; hence $PA'\neq PA$. 

(b) This now follows from the dual of (a) and Theorem \ref{THM. tangents and line-conic}(b). 
\end{proof}

\section{{\large{}Poles and polars} \label{SECTION 11. Poles and polars}}

\noindent For this section only, we adopt an additional axiom, asserting
statement (a) of Theorem \ref{THM. tangents and line-conic}. Under
Axiom P we are enabled to use Theorem \ref{Thm. 2 secants} to construct,
through any point of the plane, a secant to any conic. It remains
an open problem to determine whether this axiom may be derived from
the others. \\

\noindent \textbf{Axiom P.\label{Axiom  P}} The tangents at any three
distinct points of a point-conic are nonconcurrent; the points of
contact of any three distinct lines of a line-conic are noncollinear.
\\

The traditional method for defining a polar using a quadrangle, considering
separately points either on or not on a conic, is precluded, since
we cannot always decide, constructively, which case applies to a given
point. 
\begin{thm}
\label{Thm. Construct polar}\emph{Construction of a polar.} Let $\kappa$
be a conic, and let $P$ any point on the plane. Through the point
$P$, construct a secant $q$ of $\kappa$, using \emph{Axiom} \emph{P}
and Theorem \ref{Thm. 2 secants}. Denote the intersections of $q$
with $\kappa$ by $Q_{1}$ and $Q_{2}$, let the tangents at these
points be denoted $q_{1}$ and $q_{2}$, and set $Q=q_{1}\cdot q_{2}$.
Set $Q'=h(Q_{1},Q_{2};P)$, the harmonic conjugate of $P$ with respect
to the points $Q_{1},Q_{2}$. Then the line $p=QQ'$ is independent
of the choice of the secant $q$. 
\end{thm}

\begin{proof}
Since $Q_{1}\neq Q_{2}$, by cotransitivity we may assume that $P\neq Q_{2}$.
From Proposition \ref{Prop. tangent exists}(b) it follows that $Q_{1}\notin q_{2}$;
thus $Q_{1}Q_{2}\neq q_{2}$. From Lemma \ref{Lm. C not A. C not B.}
we have $Q'\neq Q_{2}=Q_{1}Q_{2}\cdot q_{2}$, so by Axiom C7 it follows
that $Q'\notin q_{2}$; thus $Q'\neq Q$, and the line $p$ is properly
defined. 

Now let $r$ be any secant of $\kappa$, through $P$, with $R_{1},R_{2},r_{1},r_{2},R,R'$
defined similarly, and set $s=RR'$. We must show that $s=p$.

(a) \emph{Special case; $P\in\kappa$. }In this case, $Q_{1}=P$,
so $q_{1}=t$, the tangent at $P$; thus $Q\in t$. Also, by Lemma
\ref{Lm. HC. C =00003D A}, we have $Q'=P$, so $Q'\in t$. Thus $p=t$;
similarly, $s=t$.

(b) \emph{Special case; $P\notin\kappa$ and $r\neq q$.} Since $Q_{1}\neq P=q\cdot r$,
it follows that $Q_{1}\notin r$, so $Q_{1}\neq R_{1}$. Similarly,
all four points $Q_{1},Q_{2},R_{1},R_{2}$ are distinct. By Theorem
\ref{Thm. 4 pts - tangents intersection point}, applied to the quadrangle
$Q_{1}Q_{2}R_{1}R_{2}$, the point $Q$ is collinear with the diagonals
$D_{1}=Q_{1}R_{2}\cdot Q_{2}R_{1}$ and $D_{2}=Q_{1}R_{1}\cdot Q_{2}R_{2}$,
so we have $Q\in D_{1}D_{2}$. 

The harmonic conjugate of $P$ with respect to $Q_{1},Q_{2}$ is given
by Definition \ref{Defn. HarConj}. Corresponding to the configuration
$C,A,B,l,R,P,Q,S$ in the definition, where $h(A,B;C)=AB\cdot RS$,
here we have the configuration $P,Q_{1},Q_{2},PR_{2},D_{2},R_{2},R_{1},D_{1}$.
By Theorem \ref{Thm. HC-indep}, harmonic conjugates are independent
of the choice of construction elements; thus $Q'=Q_{1}Q_{2}\cdot D_{1}D_{2}$,
so $p=D_{1}D_{2}$. Similarly, $s=D_{1}D_{2}$. 

(c) \emph{General case.} Suppose that $s\neq p$, and suppose further
that \emph{$P\notin\kappa$,} and \emph{$r\neq q$}. This contradicts
(b), negating the third assumption, so $r=q$, and now it is evident
that $s=p$. This contradicts the first assumption, negating the second;
thus $P\in\kappa$. This contradicts (a), negating the first assumption;
hence $s=p$. 
\end{proof}
\begin{defn}
\label{Defn. polar }Let $\kappa$ be a conic, and $P$ any point
on the plane. The line $p$ obtained in Theorem \ref{Thm. Construct polar}
is called the \emph{polar of $P$ with respect to $\kappa$}. 
\end{defn}

\begin{cor}
\label{Cor. polar Cor.1}Let $\kappa$ be a conic, $P$ any point
on the plane, and $p$ the polar of $P$. Then $p$ passes through:

(i) the harmonic conjugate of $P$ with respect to the points of intersection
of any secant of $\kappa$ that passes through $P$; 

(ii) the point of intersection of the tangents to $\kappa$ at the
points of intersection of any secant of $\kappa$ that passes through
$P$. 
\end{cor}

\begin{cor}
\label{Cor. polar from quadrangle}Let $\kappa$ be a conic, and let
$P$ be any point outside $\kappa$. Inscribe a quadrangle in $\kappa$
with $P$ as one diagonal point, using Theorem \ref{Thm. 2 secants}.
Then the polar of $P$ is the line joining the other two diagonal
points. 
\end{cor}

\begin{defn}
\label{Defn. pole}Let $\kappa$ be a conic, and $l$ any line on
the plane. A construction analogous to that of Theorem \ref{Thm. Construct polar}
results in a point $L$, called the \emph{pole of} $l$ \emph{with
respect to $\kappa$}. 
\end{defn}

The dual of Theorem \ref{Thm. 2 secants} yields a point $E$ on $l$,
which might be called a \emph{dual-secant}. Through $E$ pass two
distinct lines $e_{1},e_{2}$ of the line-conic $\lambda$. By Theorem
\ref{THM. tangents and line-conic}, $\lambda$ is the family of all
tangents to $\kappa$. Joining the points of contact $E_{1},E_{2}$
of the lines $e_{1},e_{2}$, we obtain a line $e=E_{1}E_{2}$. We
also construct a line $e'=h(e_{1},e_{2};l)$, the harmonic conjugate
of $l$, in the pencil $E^{*}$, with respect to the base lines $e_{1},e_{2}$.
The pole of $l$ is thus the point $L=e\cdot e'$. Corollaries analogous
to those above also apply to poles. 
\begin{thm}
\label{Thm. poles polars tangents 2 parts}Let $\kappa$ be a conic. 

(a) If a line $p$ is the polar of a point $P$, then $P$ is the
pole of $p$, and conversely. 

(b) If a point $P$ is on $\kappa$, then the polar of $P$ is the
tangent to $\kappa$ at $P$. 
\end{thm}

\begin{proof}
Statement (a) follows from Corollary \ref{Cor. polar Cor.1} and its
analog for poles. For (b), we note that in Theorem \ref{Thm. Construct polar}
for the construction of the polar $p$, we now have $Q_{1}=P$, by
Proposition \ref{Prop. conic 4  parts}(b)\emph{.} Thus $q_{1}=t$,
the tangent at $P$, and $Q\in t$. Also, by Lemma \ref{Lm. HC. C =00003D A},
we have $Q'=P$, so $Q'\in t$. Thus $p=QQ'=t$. 
\end{proof}
It remains an open problem to construct correlations and polarities
using the axioms adopted here, to develop the theory of conics constructively
using the von Staudt {[}Sta47{]} definition, whereby a conic is a
locus of points defined by a polarity, and to prove that von Staudt
conics are equivalent to the Steiner {[}Ste32{]} conics considered
here.

\part{{\large{}Analytic constructions} \label{PART II. Analytic}}

\noindent A projective plane $\mathbb{P}^{2}(\mathbb{R})$ is built
from subspaces of the linear space $\mathbb{R}^{3}$, using only constructive
properties of the real numbers. This model will establish the consistency
of the axiom system adopted in Part \ref{PART. Synthetic}. The properties
of $\mathbb{P}^{2}(\mathbb{R})$ have guided the choice of axioms,
taking note of Bishop's thesis, ``All mathematics should have numerical
meaning'' {[}B67, p. ix{]}. 

\section{{\large{}Real numbers}\label{SECTION 12. Real numbers}}

\noindent To clarify the methods used here, we give examples of familiar
properties of the real numbers that are constructively \textit{invalid},
and also properties that are constructively \textit{valid}. 

The following classical properties of a real number $\alpha$ are
constructively invalid: \emph{Either $\alpha<0$ or $\alpha=0$ or
$\alpha>0$,} and \emph{If $\neg(\alpha=0),$ then $\alpha\ne0$.}
Constructively invalid statements in classical metric geometry result
when the condition \emph{$P$ lies outside $l$,} written $P\notin l$,
is taken to mean that the distance $d(P,l)$ is positive. For examples
of constructively invalid statements for the metric plane $\mathbb{R}^{2}$,
we have the statement \emph{Either the point $P$ lies on the line
$l$ or $P$ lies outside $l$,} considered in Example \ref{EX. Pt and line - constructively invalid},
and the statement \emph{If $\neg(P\in l)$, then $P\notin l$. }

Bishop determined the constructive properties of the real numbers,
using Cauchy sequences of rationals, while referring to no axiom system
of formal logic, but only a presupposition of the positive integers
{[}B67, p.2{]}. A notable resulting feature is that the relation $\alpha\ne0$
does not refer to negation, but is given a strong affirmative definition;
one must construct an integer $n$ such that $1/n<|\alpha|$. Among
the resulting constructive properties of the reals are the following: 

(i) \emph{For any real number $\alpha$, if $\neg(\alpha\ne0)$, then
$\alpha=0$. }

(ii) \emph{For any real numbers $\alpha$ and} $\beta$, \emph{if
$\alpha\beta\neq0$, then $\alpha\neq0$ and $\beta\neq0$}. 

(iii) \emph{Given any real numbers }$\alpha$ and\emph{ $\beta$ with
}$\alpha$\emph{$<\beta$, for any real number $x$, either $x>\alpha$
or $x<\beta$.} 

Property (iii) serves as a constructive substitute for the \emph{Trichotomy
property} of classical analysis, which is constructively invalid.
For more details, and other constructive properties of the real number
system, see {[}B67, BB85, BV06{]}. For a constructive axiomatic study
of the reals, with applications to formal systems of computable analysis,
see {[}Bri99{]}. For axioms for the real numbers, and a construction
of the reals without using the axiom of countable choice, see {[}R08{]}.

\emph{Brouwerian counterexamples.} To determine the specific nonconstructivities
in a classical theory, and thereby to indicate feasible directions
for constructive work, Brouwerian counterexamples are used, in conjunction
with omniscience principles. A \emph{Brouwerian counterexample} is
a proof that a given statement implies an omniscience principle. In
turn, an \emph{omniscience principle} would imply solutions or significant
information for a large number of well-known unsolved problems. This
method was introduced by L. E. J. Brouwer {[}Bro08{]} to demonstrate
that use of the \emph{Law of Excluded Middle} inhibits mathematics
from attaining its full significance. 

Omniscience principles are primarily formulated in terms of binary
sequences, at times called \emph{decision sequences}; the zeros and
ones may represent the results of a search for a solution to a specific
problem, as in Example \ref{EX. Pt and line - constructively invalid}.
These omniscience principles have equivalent statements in terms of
real numbers; the following are those most often used in connection
with Brouwerian counterexamples. \\

\noindent \textbf{\textit{\emph{Limited principle of omniscience (LPO).
}}}\emph{For any real number $\alpha$, either $\alpha=0$ or $\alpha\neq0$.
}\\

\noindent \textbf{\textit{\emph{Lesser limited principle of omniscience
(LLPO).}}} \emph{For any real number $\alpha$, either $\alpha\leq0$
or $\alpha\geq0$.}\footnote{The omniscience principle LLPO was introduced by Bishop {[}B73{]}.}\emph{}\\

\noindent \textbf{\textit{\emph{Markov's principle (MP).}}} \emph{For
any real number $\alpha$, if $\neg(\alpha=0)$, then $\alpha\neq0$.
}\\

A statement is considered \emph{constructively invalid} if it implies
an omniscience principle. The statement considered in Example \ref{EX. Pt and line - constructively invalid},
\emph{Either} \emph{$P\in l$, or $P\notin l$,} implies LPO; thus
it is constructively invalid. For more information concerning Brouwerian
counterexamples, and other omniscience principles, see {[}B67, BB85,
M83, M88, M89, R02{]}. 

\section{{\large{}The model $\mathbb{P}^{2}(\mathbb{R})$ in Euclidean space
\label{SECTION 13. P2(R) model} }}

\noindent The model will be built following well-known classical methods,
adding constructive refinements to the definitions and proofs. 
\begin{defn}
\label{Definition.P2R-model}The plane $\mathbb{P}^{2}(\mathbb{R})$
consists of a family $\mathscr{P_{\mathrm{2}}}$ of points, and a
family $\mathscr{L_{\mathrm{2}}}$ of lines. 

\noindent $\bullet$ A \emph{point} $P$ in $\mathscr{P_{\mathrm{2}}}$
is a subspace of the linear space $\mathbb{R}^{3}$, of dimension
1. When the vector $p=(p_{1},p_{2},p_{3})$ spans $P$, we write $P=\langle p\rangle=\langle p_{1},p_{2},p_{3}\rangle$. 

\noindent $\bullet$ A \emph{line} $\lambda$ in $\mathscr{L_{\mathrm{2}}}$
is a subspace of $\mathbb{R}^{3}$, of dimension 2. When the vectors
$u,v$ span $\lambda$, and $l=u\times v$, we write $\lambda=[l]=[l_{1},l_{2},l_{3}]$,
and $\lambda=UV$, where $U=\langle u\rangle$ and $V=\langle v\rangle$. 

\noindent $\bullet$ Points $P=\left\langle p\right\rangle $ and
$Q=\langle q\rangle$ are \emph{equal}, written $P=Q$, if $p\times q=0$;
they are \emph{distinct}, written $P\neq Q$, if $p\times q\neq0$. 

\noindent $\bullet$ Lines $\lambda=[l]$ and $\mu=[m]$ are \emph{equal},
written $\lambda=\mu$, if $l\times m=0$; they are \emph{distinct},
written $\lambda\ne\mu$, if $l\times m\neq0$. 

\noindent $\bullet$ Incidence relation. Let $P=\left\langle p\right\rangle $
be a point and $\lambda=[l]$ a line. We say that $P$ \emph{lies
on} $\lambda$, and that $\lambda$ \emph{passes through} $P$, written
$P\in\lambda$, if $p\cdot l=0$. 

\noindent $\bullet$ Outside relation. For any point $P$ and any
line $\lambda$, we say that \emph{$P$} \emph{lies outside the line
}$\lambda$, and that $\lambda$ \emph{avoids the point }$P$, written
$P\notin\lambda$, if $P\neq Q$ for all points $Q$ on $\lambda$. 
\end{defn}

\noindent \emph{Notes for Definition \ref{Definition.P2R-model}. }

1. The definitions are independent of the choice of vectors spanning
the respective subspaces. That they are in accord with Definitions
\ref{Defn. Plane} and \ref{Defn. outside} for a projective plane
is evident, except for cotransitivity of the inequality relations,
which will be verified in Theorem \ref{Thm. Cotrans}. 

2. The constructive properties of the real numbers will carry over
to vectors in $\mathbb{R}^{3}$. For example, $v\neq0$ means that
at least one of the components of the vector $v$ is constructively
nonzero. 

3. The equality, inequality, incidence, and outside relations are
invariant under a change of basis. 

4. We avoid interpreting the conditions $P\in\lambda$ and $P\notin\lambda$
using the relations of set-membership and set-inclusion in the classical
sense. Theorem \ref{THM. P out ln vs dot product} and Example \ref{Examples.}
will confirm that the primary relation $P\notin\lambda$, \emph{point
outside a line,} is constructively stronger than the condition $\neg(P\in\lambda)$. 

5. The triad $\langle p_{1},p_{2},p_{3}\rangle$ is traditionally
referred to as \emph{homogeneous coördinates} for the point $P$,
and the triad $[l_{1},l_{2},l_{3}]$ as \emph{line coördinates} for
the line $\lambda$. 
\begin{thm}
\label{THM. P out ln vs dot product} Let $P=\langle p\rangle$ be
a point of $\mathbb{P}^{2}(\mathbb{R})$, and $\lambda=[l]$ a line.
Then $P\notin\lambda$ if and only if $p\cdot l\ne0$. 
\end{thm}

\begin{proof}
By a change of basis, we may assume that $l=(0,0,1)$. 

First let $P\notin\lambda$. In the case $p_{3}\neq0$ we have $p\cdot l\neq0$
directly. In the other two cases, we may set $Q=\langle q\rangle=\langle p_{1},p_{2},0\rangle$;
then $Q\in\lambda$, so $P\neq Q$. Thus $p\times q\neq0$; i.e.,
$(-p_{3}p_{2},p_{3}p_{1},0)\neq0$, and hence $p\cdot l=p_{3}\neq0$. 

For the converse, let $p\cdot l\ne0$; thus $p_{3}\neq0$. For any
point $Q=\langle q\rangle$ on $\lambda$, we have $q_{3}=0$, so
either $q_{1}\neq0$ or $q_{2}\neq0$. It follows that $p\times q=(-p_{3}q_{2},p_{3}q_{1},p_{1}q_{2}-p_{2}q_{1})\neq0$,
and hence $P\neq Q$. 
\end{proof}
\begin{cor}
Let $P=\langle p\rangle$, $Q=\langle q\rangle$, $R=\langle r\rangle$
be points of $\mathbb{P}^{2}(\mathbb{R})$, with $Q\neq R$. Then
$P\notin QR$ if and only if the vectors \textup{$p,q,r$ are independent. }
\end{cor}

\begin{cor}
\emph{\label{Cor. C6 for P2R}}Let \emph{$P=\langle p\rangle$} be
any point of $\mathbb{P}^{2}(\mathbb{R})$, and $\lambda=[u\times v]$
any line. The following conditions are equivalent: 

\emph{(a)} \emph{$\neg(P\notin\lambda)$.}

\emph{(b)} \emph{$P\in\lambda$.}

\emph{(c)} The vector $p$ is in the span of the vectors $u,v$.
\end{cor}

\begin{cor}
\label{Cor. formula for point on line}Let $Q=\langle q\rangle$,
and $R=\langle r\rangle$ be points on $\mathbb{P}^{2}(\mathbb{R})$,
with $Q\neq R$, and let $P$ be a point on $QR$. If $P\neq R$,
then there exists a unique real number $\alpha$ such that $P=\langle q+\alpha r\rangle$. 
\end{cor}

\begin{cor}
Definition \ref{Definition.P2R-model}, for the plane $\mathbb{P}^{2}(\mathbb{R})$,
is self-dual.
\end{cor}

\begin{example}
\label{Examples.}For the plane $\mathbb{P}^{2}(\mathbb{R})$, the
following statements are constructively invalid.

(a)\emph{ If $P$ and $Q$ are any points, then either $P=Q$ or $P\ne Q$}. 

(b)\emph{ If }$\lambda$\emph{ and $\mu$ are any lines, then either
$\lambda=\mu$ or $\lambda\ne\mu$}. 

(c)\emph{ If $P$ is any point, and $\lambda$ any line, then either
$P\in\lambda$ or $P\notin\lambda$.}

(d) \emph{If }$\lambda$\emph{ is any line, and $P$ is a point such
that $\neg(P\in\lambda)$, then $P\notin\lambda$}. 
\end{example}

\begin{proof}
Let \emph{$\alpha$} be any real number. For (a), set $P=\langle p\rangle=\langle0,0,1\rangle$
and $Q=\langle q\rangle=\langle\alpha,0,1\rangle$. Then $p\times q=(0,\alpha,0)$,
so the statement implies LPO. A similar counterexample serves for
(b). For (c), set $P=\langle0,\alpha,0\rangle$ and $\lambda=[0,1,0]$;
the statement implies LPO. For (d), assume also that $\neg(\alpha=0)$,
with $P$ and $\lambda$ as in the example for (c); now the statement
implies MP. 
\end{proof}
\begin{thm}
\noindent \emph{\label{Thm. Cotrans}Cotransitivity.} If $P$ and
$Q$ are points on the plane $\mathbb{P}^{2}(\mathbb{R})$, with $P\neq Q$,
then for any point $R$, either $R\neq P$, or $R\neq Q$. 
\end{thm}

\begin{proof}
\noindent By a change of basis, it suffices to consider the situation
in which $P=\langle p\rangle=\langle1,0,0\rangle$, $Q=\langle q\rangle=\langle0,1,0\rangle$,
and $R=\langle r\rangle=\langle r_{1},r_{2},r_{3}\rangle$; then $r\times p=(0,r_{3},-r_{2})$
and $r\times q=(-r_{3},0,r_{1})$. In the case $r_{1}\neq0$, we have
$r\times q\neq0$, so $R\neq Q$. In the other two cases we obtain
$R\neq P$. 
\end{proof}
\begin{lem}
\label{Lm. LinXf for projtn}On the plane $\mathbb{P}^{2}(\mathbb{R})$,
for any projection $\rho:\overline{\lambda}\rightarrow\overline{\mu}$
of a range of points $\overline{\lambda}$ onto a range $\overline{\mu}$,
there exists a non-singular linear transformation $\tau$ of $\mathbb{R}^{3}$
that induces $\rho$; i.e., $X^{\rho}=\langle\tau x\rangle$, for
all points $X=\langle x\rangle$ in the range $\overline{\lambda}$. 
\end{lem}

\begin{proof}
We adapt the proof found in {[}Art57, p. 94{]}. Select vectors $m$
and $t$ so that $\mu=[m]$ and the center of $\rho$ is $T=\langle t\rangle$.
Select distinct points $U_{1}=\langle u_{1}\rangle$ and $U_{2}=\langle u_{2}\rangle$
in $\overline{\lambda}$, select vectors $v_{i}$ such that $U_{i}^{\rho}=\langle v_{i}\rangle$,
and construct a non-singular linear transformation $\tau:\mathbb{R}^{3}\rightarrow\mathbb{R}^{3}$
such that $\tau u_{i}=v_{i}$. Since $\langle v_{1}\rangle\in TU_{1}$,
and $\langle v_{1}\rangle\neq T$, it follows from Corollary \ref{Cor. formula for point on line}
that $\langle v_{1}\rangle=\langle u_{1}+\alpha t\rangle$ for some
scalar $\alpha$, and thus $v_{1}=\beta u_{1}+t_{1}$ for a nonzero
scalar $\beta$ and a vector $t_{1}$ in $T$; we may assume that
$\beta=1$. Similarly, we have $v_{2}=u_{2}+t_{2}$, for some $t_{2}\in T$. 

Let $X=\langle x\rangle$ be any point of $\overline{\lambda}$, with
$x=\alpha_{1}u_{1}+\alpha_{2}u_{2}$. Then $\tau x=\alpha_{1}v_{1}+\alpha_{2}v_{2}$,
so $\langle\tau x\rangle\in\mu$. Also, $\tau x=\alpha_{1}u_{1}+\alpha_{2}u_{2}+t_{3}=\:x+t_{3}$,
where $t_{3}\in T$. Now set $y=t\times x$, so $TX=[y]$. Then $y=t\times(\tau x-t_{3})=t\times\tau x$,
and $\tau x\cdot y=\tau x\cdot t\times\tau x=0$, so $\langle\tau x\rangle\in TX$.
Hence $\langle\tau x\rangle=TX\cdot\mu=X^{\rho}$. 
\end{proof}
\begin{thm}
For any projectivity $\pi$ of the plane $\mathbb{P}^{2}(\mathbb{R})$,
there exists a non-singular linear transformation $\tau$ of $\mathbb{R}^{3}$
that induces $\pi$. 
\end{thm}

\begin{cor}
\label{Cor. Axiom T}Any projectivity of a range or pencil of $\mathbb{P}^{2}(\mathbb{R})$,
onto itself, with three distinct fixed elements, is the identity. 
\end{cor}

\section{{\large{}Axioms verified for the plane $\mathbb{P}^{2}(\mathbb{R})$\label{SECTION 14. Axioms P2R}}}

This verification will establish the consistency of the axiom system
adopted in Part \ref{PART. Synthetic}. The following example shows
that for Axiom C3, \emph{Distinct lines have a common point, }the
condition of distinctness is essential.
\begin{example}
\noindent \label{EX. No Com Pt}On the plane $\mathbb{P}^{2}(\mathbb{R})$,
the following statements are constructively invalid. 

(a) \emph{Given any points $P$ and $Q$, there exists a line that
passes through both points.}

(b) \emph{Given any lines $\lambda$ and $\mu$, there exists a point
that lies on both lines.}
\end{example}

\begin{proof}
It will suffice to consider the second statement. For a Brouwerian
counterexample, let $\alpha$ be any real number, and set $\alpha^{+}=\max\{\alpha,0\}$
and $\alpha^{-}=\max\{-\alpha,0\}$. Define lines $\lambda=[\alpha^{+},0,1]$
and $\mu=[0,\alpha^{-},1]$. By hypothesis, we have a point $R=\langle r\rangle=\langle r_{1},r_{2},r_{3}\rangle$
that lies on both lines. Thus $\alpha^{+}r_{1}+r_{3}=0$, and $\alpha^{-}r_{2}+r_{3}=0$.
If $r_{3}\neq0$, then we have both $\alpha^{+}\neq0$ and $\alpha^{-}\neq0$,
an absurdity; thus $r_{3}=0$. This leaves two cases. If $r_{1}\neq0$,
then $\alpha^{+}=0$, so $\alpha\leq0$, while if $r_{2}\neq0$, then
$\alpha^{-}=0$, so $\alpha\geq0$. Hence LLPO results. 
\end{proof}
\noindent \emph{Notes for Example \ref{EX. No Com Pt}. }

1. The counterexample for the first statement, the dual of that given
for the second, is easier to visualize. On $\mathbb{R}^{2}$, thought
of as a portion of $\mathbb{P}^{2}(\mathbb{R})$, consider two finite
points which are extremely near or at the origin: $P$ on the $x$-axis,
and $Q$ on the $y$-axis. If $P$ is very slightly off the origin,
and $Q$ is at the origin, then the $x$-axis is the required line
$\lambda$. But in the opposite situation we would need the $y$-axis.
Such a large jump in the output, resulting from a miniscule variation
of the input, would be a severe discontinuity in a proposed constructive
routine, and is a very strong indication that a solution would be
constructively invalid. 

2. \emph{Note on the Heyting extension }{[}H59{]}\emph{.} The above
example is essentially the same as that used in {[}M13a{]} to show
that in the Heyting extension the \emph{strong common point property}
(i.e., for all lines, not only distinct lines) is constructively invalid.
In {[}M13a, Note, p. 113{]} it was claimed that a constructive projective
plane ought to have the strong property. However, it is now seen that
various versions of a constructive real projective plane are possible.
The question of the common point property for the Heyting extension
remains an interesting open problem. For comments concerning this
issue in the classical literature, see {[}Pic75, Section 1.2{]}. 

3. \emph{Note on the projective extension of} {[}M14{]}. The strong
common point property was obtained for this plane, but the \emph{cotransitivity}
property was found to be constructively invalid {[}M14, pp. 704-5{]}.
The results of the various studies tend to indicate the incompatibility
of the two properties, strong common point and cotransitivity, in
any constructive projective plane; making this idea precise is an
open problem. 
\begin{thm}
\label{Thm. P2R. Axioms valid.}Axiom Group C, and Axioms F, D, E,
T, are valid on $\mathbb{P}^{2}(\mathbb{R})$. 
\end{thm}

\begin{proof}
Axioms C1 and C4 of Section \ref{SECTION 2. Axioms}, and Axiom E
of Section \ref{SECTION 5. Projectivities}, are evident. Given distinct
points $P=\langle p\rangle$ and $Q=\langle q\rangle$, set $l=p\times q$,
and $\lambda=[l]$. Then $p\cdot l=0$, so $P\in\lambda$, and similarly
$Q\in\lambda$. Similarly, if lines $\lambda=[l]$ and $\mu=[m]$
are distinct, then the point $P=\langle l\times m\rangle$ lies on
both lines. Thus Axioms C2 and C3 are verified. Axiom C6 was verified
in Corollary \ref{Cor. C6 for P2R}

For Axiom C5, consider lines $\lambda=[l]$ and $\mu=[m]=[u\times v]$,
and let $P=\langle p\rangle$ be a point on $\lambda$ that is outside
$\mu$. Thus $p\cdot u\times v\neq0$, and the vectors $p,u,v$ are
independent. By a change of basis, we may assume that $u,v,p$ is
the standard basis; thus $m=e_{3}$. Since $P\in\lambda$, we have
$p\cdot l=0$, so $l_{3}=0$. Now $l\times m=(l_{2},-l_{1},0)\neq0$,
and hence $\lambda\neq\mu$. 

The converse to Axiom C5 is also valid:\emph{ If the lines }$\lambda$\emph{
and }$\mu$\emph{ are distinct, then there exists a point }$P\in\lambda$\emph{
such that $P\notin\mu$.} To prove this, let $\lambda=[l]$ and $\mu=[m]$
be distinct lines. By a change of basis, we may assume that $l=e_{3}$.
Since $l\times m\neq0$, we have $(-m_{2},m_{1},0)\neq0$; thus either
$m_{1}\neq0$ or $m_{2}\neq0$. In the first case, set $P=\langle e_{1}\rangle$.
Then $p\cdot l=0$, so $P\in\lambda$, while $p\cdot m=m_{1}\neq0$,
so $P\notin\mu$. The second case is similar. 

For Axiom C7, let $\lambda$ and $\mu$ be distinct lines, and let
$P\neq Q=\lambda\cdot\mu$. Select a point $R\in\lambda$ so that
$R\neq Q$; thus $\lambda=QR.$ By the converse to Axiom C5, verified
in the preceding paragraph, we may select a point $S\in\mu$ such
that $S\notin\lambda$. Thus $S\neq Q$, so $\mu=QS,$ and the points
$Q,R,S$ are noncollinear. By a change of basis, we may assume that
$R=\langle e_{1}\rangle$, $S=\langle e_{2}\rangle$, $Q=\langle e_{3}\rangle$,
and $P=\langle p\rangle=\langle p_{1},p_{2},p_{3}\rangle$. Then $\lambda=[e_{2}]$,
and $\mu=[e_{1}]$. Since $P\neq Q$, we have $p\times e_{3}\neq0$,
so $(p_{2},-p_{1},0)\neq0$, and thus either $p_{2}\neq0$, or $p_{1}\neq0$.
In the first case, we have $P\notin\lambda$, while in the second
case we find that $P\notin\mu$. 

For Axiom F of Section \ref{SECTION 2. Axioms}, Fano's Axiom, by
a change of basis we may assume that the quadrangle $PQRS$ has vertices
$\langle e_{1}\rangle,\langle e_{2}\rangle,\langle e_{3}\rangle,$
and $\langle e\rangle=\langle1,1,1\rangle$. The six sides are then
$PQ=[0,0,1]$, $PR=[0,1,0]$, $PS=[0,1,-1]$, $QR=[1,0,0]$, $QS=[1,0,-1]$,
and $RS=[1,-1,0]$. The diagonal points are $D_{1}=PQ\cdot RS=\langle1,1,0\rangle$,
$D_{2}=PR\cdot QS=\langle1,0,1\rangle$, and $D_{3}=\langle d\rangle=PS\cdot QR=\langle0,1,1\rangle$.
Thus $D_{1}D_{2}=[m]=[1,-1,-1]$. Since $d\cdot m\neq0$, we have
$D_{3}\notin D_{1}D_{2}$. Thus the diagonal points are noncollinear.

For Axiom D of Section \ref{SECTION 3. Des}, Desargues's Theorem,
consider triangles $PQR$ and $P'Q'R'$, perspective from a center
$O=\langle o\rangle$. By a change of basis, we may assume that $P,Q,R=\langle e_{1}\rangle,\langle e_{2}\rangle,\langle e_{3}\rangle$;
thus $QR=[e_{1}]$, $RP=[e_{2}]$, and $PQ=[e_{3}]$. Since $O\notin PQ$,
we have $o_{3}\neq0$, and similarly for $o_{1}$ and $o_{2}$; thus
we may assume that $O=\langle1,1,1\rangle$. Now $OP=[0,1,-1]$. Since
$P'\in OP$, and $P'\neq P$, it follows from Corollary \ref{Cor. formula for point on line}
that $P'=\langle\alpha,1,1\rangle$ for some scalar $\alpha$; similarly,
$Q'=\langle1,\beta,1\rangle$, and $R'=\langle1,1,\gamma\rangle$.
It follows that $P'Q'=[1-\beta,1-\alpha,\alpha\beta-1]$, $Q'R'=[\beta\gamma-1,1-\gamma,1-\beta]$,
and $R'P'=[1-\gamma,\gamma\alpha-1,1-\alpha]$. Now the points in
question are $A=\langle a\rangle=PQ\cdot P'Q'=\langle\alpha-1,1-\beta,0\rangle$,
$B=\langle b\rangle=QR\cdot Q'R'=\langle0,\beta-1,1-\gamma\rangle$,
and $C=\langle c\rangle=RP\cdot R'P'=\langle\alpha-1,0,1-\gamma\rangle$.
Since $c\cdot a\times b$ works out to $0$, it follows that $C\in AB$.
Thus the points $A,B,C$ are collinear. To show that the line $AB$
avoids each of the six vertices, we first note that since the center
$O$ lies outside each of the six sides, we have $O\neq Q'$; it follows
that $\beta\neq1$, and similarly, $\gamma\neq1$. Thus $e_{1}\cdot a\times b=(1-\beta)(1-\gamma)\neq0,$
and we have $P\notin AB$. By symmetry of the vertices of the triangle,
$Q$ and $R$ also lie outside $AB$. By symmetry of the two triangles,
the points $P',Q',R'$ lie outside $AB$. Hence the triangles are
perspective from the axis $AB$. 

Axiom T of Section \ref{SECTION 6. Fund Thm}, the uniqueness portion
of the \emph{Fundamental Theorem}, was verified in Corollary \ref{Cor. Axiom T}. 
\end{proof}
It remains an open problem to develop the analytic theory of conics
constructively, and to determine the constructive validity of Axiom
P of Section \ref{SECTION 11. Poles and polars} in an analytic setting.
\\

\noindent \textbf{Acknowledgments.} The author is grateful to the
referee and to Dr. G. Calderón for many useful suggestions. \\

\section*{{\large{}References}{\small{} \label{REFERENCES}}\protect \\
}

\noindent {\small{}{[}Art57{]} E. Artin,}\emph{\small{} Geometric
Algebra,}{\small{} Interscience, New York, 1957. MR0082463 }{\small\par}

\noindent {\small{}{[}B65{]} E. Bishop, Book Review: }\emph{\small{}The
Foundations of Intuitionistic Mathematics, }{\small{}by S. C. Kleene
and R. E. Vesley, Bull. Amer. Math. Soc. 71:850-852, 1965. MR1566375 }{\small\par}

\noindent {\small{}{[}B67{]}E. Bishop, }\emph{\small{}Foundations
of Constructive Analysis,}{\small{} McGraw-Hill, New York, 1967. MR0221878 }{\small\par}

\noindent {\small{}{[}B73{]} E. Bishop, }\emph{\small{}Schizophrenia
in Contemporary Mathematics}{\small{}, AMS Colloquium Lectures, Missoula,
Montana, 1973. Reprinted in }\emph{\small{}Contemporary Mathematics}{\small{}
39:l-32, 1985. MR0788163}{\small\par}

\noindent {\small{}{[}B75{]} E. Bishop, }\emph{\small{}The crisis
in contemporary mathematics,}{\small{} Proceedings of the American
Academy Workshop on the Evolution of Modern Mathematics, Boston, Mass.,
1974; }\emph{\small{}Historia Math.}{\small{} 2:507\textendash 517,
1975. MR0498014 }{\small\par}

\noindent {\small{}{[}BB85{]} E. Bishop and D. Bridges, }\emph{\small{}Constructive
Analysis,}{\small{} Springer-Verlag, Berlin, 1985. MR0804042 }{\small\par}

\noindent {\small{}{[}BR87{]} D. Bridges and F. Richman, }\emph{\small{}Varieties
of Constructive Mathematics,}{\small{} Cambridge University Press,
Cambridge, UK, 1987. MR0890955 }{\small\par}

\noindent {\small{}{[}BV06{]} D. Bridges and L. Vita, }\emph{\small{}Techniques
of Constructive Analysis,}{\small{} Springer, New York, 2006. MR2253074 }{\small\par}

\noindent {\small{}{[}Bee10{]} M. Beeson, Constructive geometry, }\emph{\small{}10th
Asian Logic Conference}{\small{}, pp. 19\textendash 84, World Sci.
Publ., Hackensack, NJ, 2010. MR2798893 }{\small\par}

\noindent {\small{}{[}Bri99{]} D. Bridges, Constructive mathematics:
a foundation for computable analysis, }\emph{\small{}Theoret. Comput.
Sci.}{\small{} 219:95\textendash 109, 1999. MR1694428 }{\small\par}

\noindent {\small{}{[}Bro08{]} L. E. J. Brouwer, De onbetrouwbaarheid
der logische principes, }\emph{\small{}Tijdschrift voor Wijsbegeerte}{\small{}
2:152-158, 1908. English translation; The Unreliability of the Logical
Principles, pp. 107\textendash 111 in A. Heyting, ed., }\emph{\small{}L.
E. J. Brouwer: Collected Works 1: Philosophy and Foundations of Mathematics,}{\small{}
Elsevier, Amsterdam-New York, 1975. MR0532661 }{\small\par}

\noindent {\small{}{[}Cox55{]} H. S. M. Coxeter, }\emph{\small{}The
Real Projective Plane, 2nd ed.,}{\small{} University Press, Cambridge,
1955. MR0070189 }{\small\par}

\noindent {\small{}{[}Cre73{]} L. Cremona, }\emph{\small{}Elementi
di geometria projettiva, }{\small{}G. B. Paravia e Comp., Torino,
1873. English translation; }\emph{\small{}Elements of Projective Geometry,}{\small{}
trans. by C. Leudesdorf, Clarendon Press, Oxford, 1985. Reprint; Forgotten
Books, Hong Kong, 2012. }{\small\par}

\noindent {\small{}{[}D63{]} D. van Dalen, Extension problems in intuitionistic
plane projective geometry I, II, }\emph{\small{}Indag. Math.}{\small{}
25:349-383, 1963. MR0153567 MR0153568 }{\small\par}

\noindent {\small{}{[}D90{]} D. van Dalen, Heyting and intuitionistic
geometry, pp. 19-27 in }\emph{\small{}Mathematical logic,}{\small{}
Plenum, New York, 1990. MR1083983 }{\small\par}

\noindent {\small{}{[}D96{]} D. van Dalen, 'Outside' as a primitive
notion in constructive projective geometry, }\emph{\small{}Geom. Dedicata}{\small{}
60:107-111, 1996. MR1376483 }{\small\par}

\noindent {\small{}{[}Fan92{]} G. Fano, Sui postulati fondamentali
della geometria proiettiva in uno spazio lineare a un numero qualunque
di dimensioni, }\emph{\small{}Giornale di mat. di Battista}{\small{}
30:106\textendash 132, 1892. }{\small\par}

\noindent {\small{}{[}H28{]} A. Heyting, Zur intuitionistischen Axiomatik
der projektiven Geometrie, }\emph{\small{}Math. Ann.}{\small{} 98:491-538,
1928. MR1512416}{\small\par}

\noindent {\small{}{[}H59{]} A. Heyting, Axioms for intuitionistic
plane affine geometry, pp. 160-173 in L. Henkin, P. Suppes, A. Tarski,
eds., }\emph{\small{}The Axiomatic Method, with special reference
to geometry and physics: Proceedings of an international symposium
held at the University of California, Berkeley, December 26, 1957
- January 4, 1958,}{\small{} North-Holland, Amsterdam, 1959. MR0120154 }{\small\par}

\noindent {\small{}{[}H66{]} A. Heyting, }\emph{\small{}Intuitionism:
An Introduction,}{\small{} North-Holland, Amsterdam, 1966. MR0221911 }{\small\par}

\noindent {\small{}{[}Kli72{]} M. Kline, }\emph{\small{}Mathematical
Thought from Ancient to Modern Times}{\small{}, Oxford University
Press, 1972. MR0472307 }{\small\par}

\noindent {\small{}{[}Leh17{]} D. N. Lehmer, }\emph{\small{}An Elementary
Course in Synthetic Projective Geometry,}{\small{} Ginn, Boston, 1917. }{\small\par}

\noindent {\small{}{[}LomVes98{]} M. Lombard and R. Vesley, A common
axiom set for classical and intuitionistic plane geometry, }\emph{\small{}Ann.
Pure Appl. Logic}{\small{} 95:229-255, 1998. MR1650655 }{\small\par}

\noindent {\small{}{[}M83{]} M. Mandelkern, }\emph{\small{}Constructive
Continuity,}{\small{} Mem. Amer. Math. Soc. 42(277), 1983. MR0690901}{\small\par}

\noindent {\small{}{[}M85{]} M. Mandelkern, Constructive mathematics,
}\emph{\small{}Math. Mag. }{\small{}58:272-280, 1985. MR0810148 }{\small\par}

\noindent {\small{}{[}M88{]} M. Mandelkern, Limited omniscience and
the Bolzano-Weierstrass principle, }\emph{\small{}Bull. London Math.
Soc.}{\small{} 20:319-320, 1988. MR0940284 }{\small\par}

\noindent {\small{}{[}M89{]} M. Mandelkern, Brouwerian counterexamples,
}\emph{\small{}Math. Mag.}{\small{} 62:3-27, 1989. MR0986618 }{\small\par}

\noindent {\small{}{[}M07{]} M. Mandelkern, Constructive coördinatization
of Desarguesian planes, }\emph{\small{}Beiträge Algebra Geom.}{\small{}
48:547-589, 2007. MR2364806 }{\small\par}

\noindent {\small{}{[}M13a{]} M. Mandelkern, The common point problem
in constructive projective geometry, }\emph{\small{}Indag. Math.}{\small{}
}\emph{\small{}(N.S.)}{\small{} 24:111-114, 2013. MR2997755}{\small\par}

\noindent {\small{}{[}M13b{]} M. Mandelkern, Convexity and osculation
in normed spaces, }\emph{\small{}Rocky Mountain J. Math.}{\small{}
43:551-561, 2013. MR3077842 }{\small\par}

\noindent {\small{}{[}M14{]} M. Mandelkern, Constructive projective
extension of an incidence plane, }\emph{\small{}Trans. Amer. Math.
Soc.}{\small{} 366:691-706, 2014. MR3130314 }{\small\par}

\noindent {\small{}{[}Pam01{]} V. Pambuccian, Constructive axiomatization
of plane hyperbolic geometry, }\emph{\small{}MLQ Math. Log. Q.}{\small{}
47:475\textendash 488, 2001. MR1865767 }{\small\par}

\noindent {\small{}{[}Pam11{]} V. Pambuccian, The axiomatics of ordered
geometry; I. Ordered incidence spaces, }\emph{\small{}Expo. Math.}{\small{}
29:24\textendash 66, 2011. MR2785544}{\small\par}

\noindent {\small{}{[}Pic75{]} G. Pickert, }\emph{\small{}Projektive
Ebenen}{\small{}, }\emph{\small{}2. Aufl.,}{\small{} Springer-Verlag,
New York, Berlin, 1975. MR0370350 }{\small\par}

\noindent {\small{}{[}Pla95{]} J. von Plato, The axioms of constructive
geometry, }\emph{\small{}Ann. Pure Appl. Logic}{\small{} 76:169-200,
1995. MR1361484 }{\small\par}

\noindent {\small{}{[}Pla98{]} J. von Plato, A constructive theory
of ordered affine geometry, }\emph{\small{}Indag. Math.}{\small{}
}\emph{\small{}(N.S.)}{\small{} 9:549-562, 1998. MR1691994 }{\small\par}

\noindent {\small{}{[}Pla10{]} J. von Plato, Combinatorial analysis
of proofs in projective and affine geometry, }\emph{\small{}Ann. Pure
Appl. Logic}{\small{} 162:144-161, 2010. MR2737928 }{\small\par}

\noindent {\small{}{[}Pon22{]} J. V. Poncelet, }\emph{\small{}Traité
des Propriétés Projectives des Figures,}{\small{} Gauthier-Villars,
Paris, 1822. }{\small\par}

\noindent {\small{}{[}R82{]} F. Richman, Meaning and information in
constructive mathematics, }\emph{\small{}Amer. Math. Monthly}{\small{}
89:385-388, 1982. MR0660918 }{\small\par}

\noindent {\small{}{[}R02{]} F. Richman, Omniscience principles and
functions of bounded variation, }\emph{\small{}MLQ Math. Log. Q.}{\small{}
42:111-116, 2002. MR1874208 }{\small\par}

\noindent {\small{}{[}R08{]} F. Richman, Real numbers and other completions,
MLQ Math. Log. Q. 54:98\textendash 108, 2008. MR2387400 }{\small\par}

\noindent {\small{}{[}Sta47{]} K. G. C. von Staudt,}\emph{\small{}
Geometrie der Lage, }{\small{}Nürnberg, 1847. }{\small\par}

\noindent {\small{}{[}Ste32{]} J. Steiner, }\emph{\small{}Sytematische
Entwickelung der Abhängigkeit geometrischer Gestalten von einander,
}{\small{}1832. }{\small\par}

\noindent {\small{}{[}Sto70{]} G. Stolzenberg, Book Review: }\emph{\small{}Foundations
of Constructive Analysis,}{\small{} by E. Bishop,}\emph{\small{} Bull.
Amer. Math. Soc.}{\small{} 76:301-323, 1970. MR0257687 }{\small\par}

\noindent {\small{}{[}VY10{]} O. Veblen and J. W. Young,}\emph{\small{}
Projective Geometry, Vol. 1,}{\small{} Ginn, Boston, 1910. MR0179666 }{\small\par}

\noindent {\small{}{[}Wei07{]} C. Weibel, Survey of Non-Desarguesian
Planes, }\emph{\small{}Notices Amer. Math. Soc.}{\small{} 54:1294\textendash 1303,
2007. }{\small\par}

\noindent {\small{}{[}You30{]} J. W. Young, }\emph{\small{}Projective
Geometry,}{\small{} Math. Assoc. Amer., Open Court, Chicago, 1930.
}\\

\noindent {\small{}New Mexico State University}{\small\par}

\noindent {\small{}Las Cruces, New Mexico }{\small\par}

\noindent \emph{\small{}e-mail:}{\small{} mandelkern@zianet.com }{\small\par}

\noindent \emph{\small{}web:}{\small{} http://www.zianet.com/mandelkern
}\\

\noindent {\small{}August 12, 2014; revised, February 14, 2015. }\\

\end{document}